\documentclass[a4paper]{article}
\usepackage{amsmath,amsfonts,amssymb,amscd,graphicx,bez123}

\title{The Dynkin diagrams of Rational Double Points}
\author{Benjamin Friedrich\thanks{Supervisor: Prof. P.M.H. Wilson}}

\date{15. May 2002}

\begin{document}

% A2 - big
\newcommand{\Azweibig}
{
\setlength{\unitlength}{1mm}
\begin{picture}(15,10)(0,-5)

\put (0,0) {\circle{3}}

\put (2,0) {\line(1,0){6}}

\put (10,0) {\circle{3}}

\end{picture}
}

% E6
\newcommand{\Esechs}
{
\setlength{\unitlength}{0.6mm}
\begin{picture}(50,30)(-5,-20)

\put (0,0) {\circle{3}}

\put (-3,4) {\makebox{-2}}

\put (2,0) {\line(1,0){6}}

\put (10,0) {\circle{3}}

\put (7,4) {\makebox{-2}}

\put (12,0) {\line(1,0){6}}

\put (20,0) {\circle{3}}

\put (17,4) {\makebox{-2}}

\put (22,0) {\line(1,0){6}}

\put (30,0) {\circle{3}}

\put (27,4) {\makebox{-2}}

\put (32,0) {\line(1,0){6}}

\put (40,0) {\circle{3}}

\put (37,4) {\makebox{-2}}

\put (20,-8) {\line(0,1){6}}

\put (20,-10) {\circle{3}}

\put (17,-17) {\makebox{-2}}

\end{picture}
}

% E6 - Znum
\newcommand{\Esechsznum}
{
\setlength{\unitlength}{0.6mm}
\begin{picture}(50,30)(-5,-20)

\put (0,0) {\circle{3}}

\put (-2,4) {\makebox{1}}

\put (2,0) {\line(1,0){6}}

\put (10,0) {\circle{3}}

\put (8,4) {\makebox{2}}

\put (12,0) {\line(1,0){6}}

\put (20,0) {\circle{3}}

\put (18,4) {\makebox{3}}

\put (22,0) {\line(1,0){6}}

\put (30,0) {\circle{3}}

\put (28,4) {\makebox{2}}

\put (32,0) {\line(1,0){6}}

\put (40,0) {\circle{3}}

\put (38,4) {\makebox{1}}

\put (20,-8) {\line(0,1){6}}

\put (20,-10) {\circle{3}}

\put (18,-17) {\makebox{2}}

\end{picture}
}

% E7
\newcommand{\Esieben}
{
\setlength{\unitlength}{0.6mm}
\begin{picture}(60,30)(-5,-20)

\put (0,0) {\circle{3}}

\put (-3,4) {\makebox{-2}}

\put (2,0) {\line(1,0){6}}

\put (10,0) {\circle{3}}

\put (7,4) {\makebox{-2}}

\put (12,0) {\line(1,0){6}}

\put (20,0) {\circle{3}}

\put (17,4) {\makebox{-2}}

\put (22,0) {\line(1,0){6}}

\put (30,0) {\circle{3}}

\put (27,4) {\makebox{-2}}

\put (32,0) {\line(1,0){6}}

\put (40,0) {\circle{3}}

\put (37,4) {\makebox{-2}}

\put (42,0) {\line(1,0){6}}

\put (50,0) {\circle{3}}

\put (47,4) {\makebox{-2}}

\put (20,-8) {\line(0,1){6}}

\put (20,-10) {\circle{3}}

\put (17,-17) {\makebox{-2}}

\end{picture}
}

% E7 - Znum
\newcommand{\Esiebenznum}
{
\setlength{\unitlength}{0.6mm}
\begin{picture}(60,30)(-5,-20)

\put (0,0) {\circle{3}}

\put (-2,4) {\makebox{2}}

\put (2,0) {\line(1,0){6}}

\put (10,0) {\circle{3}}

\put (8,4) {\makebox{3}}

\put (12,0) {\line(1,0){6}}

\put (20,0) {\circle{3}}

\put (18,4) {\makebox{4}}

\put (22,0) {\line(1,0){6}}

\put (30,0) {\circle{3}}

\put (28,4) {\makebox{3}}

\put (32,0) {\line(1,0){6}}

\put (40,0) {\circle{3}}

\put (38,4) {\makebox{2}}

\put (42,0) {\line(1,0){6}}

\put (50,0) {\circle{3}}

\put (48,4) {\makebox{1}}

\put (20,-8) {\line(0,1){6}}

\put (20,-10) {\circle{3}}

\put (18,-17) {\makebox{2}}

\end{picture}
}

% E8
\newcommand{\Eacht}
{
\setlength{\unitlength}{0.6mm}
\begin{picture}(70,30)(-5,-20)

\put (0,0) {\circle{3}}

\put (-3,4) {\makebox{-2}}

\put (2,0) {\line(1,0){6}}

\put (10,0) {\circle{3}}

\put (7,4) {\makebox{-2}}

\put (12,0) {\line(1,0){6}}

\put (20,0) {\circle{3}}

\put (17,4) {\makebox{-2}}

\put (22,0) {\line(1,0){6}}

\put (30,0) {\circle{3}}

\put (27,4) {\makebox{-2}}

\put (32,0) {\line(1,0){6}}

\put (40,0) {\circle{3}}

\put (37,4) {\makebox{-2}}

\put (42,0) {\line(1,0){6}}

\put (50,0) {\circle{3}}

\put (47,4) {\makebox{-2}}

\put (52,0) {\line(1,0){6}}

\put (60,0) {\circle{3}}

\put (57,4) {\makebox{-2}}

\put (20,-8) {\line(0,1){6}}

\put (20,-10) {\circle{3}}

\put (17,-17) {\makebox{-2}}

\end{picture}
}

% E8 - Znum
\newcommand{\Eachtznum}
{
\setlength{\unitlength}{0.6mm}
\begin{picture}(70,30)(-5,-20)

\put (0,0) {\circle{3}}

\put (-2,4) {\makebox{2}}

\put (2,0) {\line(1,0){6}}

\put (10,0) {\circle{3}}

\put (8,4) {\makebox{4}}

\put (12,0) {\line(1,0){6}}

\put (20,0) {\circle{3}}

\put (18,4) {\makebox{6}}

\put (22,0) {\line(1,0){6}}

\put (30,0) {\circle{3}}

\put (28,4) {\makebox{5}}

\put (32,0) {\line(1,0){6}}

\put (40,0) {\circle{3}}

\put (38,4) {\makebox{4}}

\put (42,0) {\line(1,0){6}}

\put (50,0) {\circle{3}}

\put (48,4) {\makebox{3}}

\put (52,0) {\line(1,0){6}}

\put (60,0) {\circle{3}}

\put (58,4) {\makebox{2}}

\put (20,-8) {\line(0,1){6}}

\put (20,-10) {\circle{3}}

\put (18,-17) {\makebox{3}}

\end{picture}
}

% Dn
\newcommand{\Dn}
{
\setlength{\unitlength}{0.6mm}
\begin{picture}(70,30)(-5,-20)

\put (0,0) {\circle{3}}

\put (-3,4) {\makebox{-2}}

\put (2,0) {\line(1,0){6}}

\put (10,0) {\circle{3}}

\put (7,4) {\makebox{-2}}

\put (12,0) {\line(1,0){6}}

\put (20,0) {\circle{3}}

\put (17,4) {\makebox{-2}}

\put (22,0) {\line(1,0){6}}

\put (32,0) {\circle*{2}}
\put (35,0) {\circle*{2}}
\put (38,0) {\circle*{2}}

\put (42,0) {\line(1,0){6}}

\put (50,0) {\circle{3}}

\put (47,4) {\makebox{-2}}

\put (52,0) {\line(1,0){6}}

\put (60,0) {\circle{3}}

\put (57,4) {\makebox{-2}}

\put (10,-8) {\line(0,1){6}}

\put (10,-10) {\circle{3}}

\put (8,-17) {\makebox{-2}}

\end{picture}
}

% Dn - Znum
\newcommand{\Dnznum}
{
\setlength{\unitlength}{0.6mm}
\begin{picture}(70,30)(-5,-20)

\put (0,0) {\circle{3}}

\put (-2,4) {\makebox{1}}

\put (2,0) {\line(1,0){6}}

\put (10,0) {\circle{3}}

\put (8,4) {\makebox{2}}

\put (12,0) {\line(1,0){6}}

\put (20,0) {\circle{3}}

\put (18,4) {\makebox{2}}

\put (22,0) {\line(1,0){6}}

\put (32,0) {\circle*{2}}
\put (35,0) {\circle*{2}}
\put (38,0) {\circle*{2}}

\put (42,0) {\line(1,0){6}}

\put (50,0) {\circle{3}}

\put (48,4) {\makebox{2}}

\put (52,0) {\line(1,0){6}}

\put (60,0) {\circle{3}}

\put (58,4) {\makebox{1}}

\put (10,-8) {\line(0,1){6}}

\put (10,-10) {\circle{3}}

\put (8,-17) {\makebox{1}}

\end{picture}
}

% An
\newcommand{\An}
{
\setlength{\unitlength}{0.6mm}
\begin{picture}(70,15)(-5,-5)

\put (0,0) {\circle{3}}

\put (-3,4) {\makebox{-2}}

\put (2,0) {\line(1,0){6}}

\put (10,0) {\circle{3}}

\put (7,4) {\makebox{-2}}

\put (12,0) {\line(1,0){6}}

\put (20,0) {\circle{3}}

\put (17,4) {\makebox{-2}}

\put (22,0) {\line(1,0){6}}

\put (32,0) {\circle*{2}}
\put (35,0) {\circle*{2}}
\put (38,0) {\circle*{2}}

\put (42,0) {\line(1,0){6}}

\put (50,0) {\circle{3}}

\put (47,4) {\makebox{-2}}

\put (52,0) {\line(1,0){6}}

\put (60,0) {\circle{3}}

\put (57,4) {\makebox{-2}}

\end{picture}
}

% An - Znum
\newcommand{\Anznum}
{
\setlength{\unitlength}{0.6mm}
\begin{picture}(70,15)(-5,-5)

\put (0,0) {\circle{3}}

\put (-2,4) {\makebox{1}}

\put (2,0) {\line(1,0){6}}

\put (10,0) {\circle{3}}

\put (8,4) {\makebox{1}}

\put (12,0) {\line(1,0){6}}

\put (20,0) {\circle{3}}

\put (18,4) {\makebox{1}}

\put (22,0) {\line(1,0){6}}

\put (32,0) {\circle*{2}}
\put (35,0) {\circle*{2}}
\put (38,0) {\circle*{2}}

\put (42,0) {\line(1,0){6}}

\put (50,0) {\circle{3}}

\put (48,4) {\makebox{1}}

\put (52,0) {\line(1,0){6}}

\put (60,0) {\circle{3}}

\put (58,4) {\makebox{1}}

\end{picture}
}

% cross
\newcommand{\cross}
{
\setlength{\unitlength}{1mm}
\begin{picture}(80,30)(-15,-15)

\put (0,0) {\circle{3}}

\put (-3,4) {\makebox{$e_1$}}

\put (2,0) {\line(1,0){6}}

\put (10,0) {\circle{3}}

\put (7,4) {\makebox{$e_2$}}

\put (12,0) {\line(1,0){6}}

\put (22,0) {\circle*{1}}
\put (25,0) {\circle*{1}}
\put (28,0) {\circle*{1}}

\put (32,0) {\line(1,0){6}}

\put (40,0) {\circle{3}}

\put (37,4) {\makebox{$e_{n-1}$}}

\put (42,0) {\line(1,0){6}}

\put (50,0) {\circle{3}}

\put (47,4) {\makebox{$e_n$}}

\put (52,2) {\line(1,1){4}}

\put (58,8) {\circle{3}}

\put (60,7) {\makebox{$f_3$}}

\put (52,-2) {\line(1,-1){4}}

\put (58,-8) {\circle{3}}

\put (60,-9) {\makebox{$f_4$}}

\put (-6,6) {\line(1,-1){4}}

\put (-8,8) {\circle{3}}

\put (-14,7) {\makebox{$f_1$}}

\put (-6,-6) {\line(1,1){4}}

\put (-8,-8) {\circle{3}}

\put (-14,-9) {\makebox{$f_2$}}

\end{picture}
}

% T-tree
\newcommand{\Ttree}
{
\setlength{\unitlength}{0.7mm}
\begin{picture}(120,60)(-5,-50)

\put (0,0) {\circle{3}}

\put (2,0) {\line(1,0){6}}

\put (10,0) {\circle{3}}

\put (12,0) {\line(1,0){6}}

\put (22,0) {\circle*{1}}
\put (25,0) {\circle*{1}}
\put (28,0) {\circle*{1}}

\put (32,0) {\line(1,0){6}}

\put (40,0) {\circle{3}}

\put (42,0) {\line(1,0){6}}

\put (50,0) {\circle{3}}

\put (52,0) {\line(1,0){6}}

\put (60,0) {\circle{3}}

\put (62,0) {\line(1,0){6}}

\put (72,0) {\circle*{1}}
\put (75,0) {\circle*{1}}
\put (78,0) {\circle*{1}}

\put (82,0) {\line(1,0){6}}

\put (90,0) {\circle{3}}

\put (92,0) {\line(1,0){6}}

\put (100,0) {\circle{3}}

\put (50,-8) {\line(0,1){6}}

\put (50,-10) {\circle{3}}

\put (50,-18) {\line(0,1){6}}

\put (50,-22) {\circle*{1}}
\put (50,-25) {\circle*{1}}
\put (50,-28) {\circle*{1}}

\put (50,-38) {\line(0,1){6}}

\put (50,-40) {\circle{3}}

\end{picture}
}

% cyclus
\newcommand{\cyclus}
{
\setlength{\unitlength}{1mm}
\begin{picture}(30,25)(-17,-5)

\put (0,0) {\circle{3}}

\put (-1,-5) {\makebox{$e_1$}}

\put (2,2) {\line(1,1){4}}

\put (8,8) {\circle{3}}

\put (10,7) {\makebox{$e_2$}}

\put (8,10) {\line(0,1){6}}

\put (8,18) {\circle{3}}

\put (10,17) {\makebox{$e_3$}}

\put (-3,18) {\circle*{1}}
\put (0,18) {\circle*{1}}
\put (3,18) {\circle*{1}}

\put (-6,6) {\line(1,-1){4}}

\put (-8,8) {\circle{3}}

\put (-14,7) {\makebox{$e_n$}}

\put (-8,10) {\line(0,1){6}}

\put (-8,18) {\circle{3}}

\put (-17,17) {\makebox{$e_{n-1}$}}

\end{picture}
}

% Hirzebruch - Jung
\newcommand{\HiJu}
{
\setlength{\unitlength}{1mm}
\begin{picture}(60,20)(-5,-5)

\put (0,0) {\circle{3}}

\put (-4,4) {\makebox{$-b_1$}}

\put (2,0) {\line(1,0){6}}

\put (10,0) {\circle{3}}

\put (6,4) {\makebox{$-b_2$}}

\put (12,0) {\line(1,0){6}}

\put (20,0) {\circle{3}}

\put (16,4) {\makebox{$-b_3$}}

\put (22,0) {\line(1,0){6}}

\put (32,0) {\circle*{1}}
\put (35,0) {\circle*{1}}
\put (38,0) {\circle*{1}}

\put (42,0) {\line(1,0){6}}

\put (50,0) {\circle{3}}

\put (46,4) {\makebox{$-b_k$}}

\end{picture}
}

% E8 labelled
\newcommand{\Eachtlab}
{
\setlength{\unitlength}{1mm}
\begin{picture}(70,30)(-5,-20)

\put (0,0) {\circle{5}}

\put (-2,4) {\makebox{-2}}

\put (-2,-1) {\makebox{$e_1$}}

\put (3,0) {\line(1,0){4}}

\put (10,0) {\circle{5}}

\put (8,4) {\makebox{-2}}

\put (8,-1) {\makebox{$e_2$}}

\put (13,0) {\line(1,0){4}}

\put (20,0) {\circle{5}}

\put (18,4) {\makebox{-2}}

\put (18,-1) {\makebox{$e_3$}}

\put (23,0) {\line(1,0){4}}

\put (30,0) {\circle{5}}

\put (28,4) {\makebox{-2}}

\put (28,-1) {\makebox{$e_5$}}

\put (33,0) {\line(1,0){4}}

\put (40,0) {\circle{5}}

\put (38,4) {\makebox{-2}}

\put (38,-1) {\makebox{$e_6$}}

\put (43,0) {\line(1,0){4}}

\put (50,0) {\circle{5}}

\put (48,4) {\makebox{-2}}

\put (48,-1) {\makebox{$e_7$}}

\put (53,0) {\line(1,0){4}}

\put (60,0) {\circle{5}}

\put (58,4) {\makebox{-2}}

\put (58,-1) {\makebox{$e_8$}}

\put (20,-7) {\line(0,1){4}}

\put (20,-10) {\circle{5}}

\put (18,-17) {\makebox{-2}}

\put (18,-11) {\makebox{$e_4$}}

\end{picture}
}

% E6 - plain
\newcommand{\Esechsplain}
{
\setlength{\unitlength}{0.6mm}
\begin{picture}(50,20)(-5,-15)

\put (0,0) {\circle{3}}

\put (2,0) {\line(1,0){6}}

\put (10,0) {\circle{3}}

\put (12,0) {\line(1,0){6}}

\put (20,0) {\circle{3}}

\put (22,0) {\line(1,0){6}}

\put (30,0) {\circle{3}}

\put (32,0) {\line(1,0){6}}

\put (40,0) {\circle{3}}

\put (20,-8) {\line(0,1){6}}

\put (20,-10) {\circle{3}}

\end{picture}
}

% E7 - plain
\newcommand{\Esiebenplain}
{
\setlength{\unitlength}{0.6mm}
\begin{picture}(60,20)(-5,-15)

\put (0,0) {\circle{3}}

\put (2,0) {\line(1,0){6}}

\put (10,0) {\circle{3}}

\put (12,0) {\line(1,0){6}}

\put (20,0) {\circle{3}}

\put (22,0) {\line(1,0){6}}

\put (30,0) {\circle{3}}

\put (32,0) {\line(1,0){6}}

\put (40,0) {\circle{3}}

\put (42,0) {\line(1,0){6}}

\put (50,0) {\circle{3}}

\put (20,-8) {\line(0,1){6}}

\put (20,-10) {\circle{3}}

\end{picture}
}

% E8 - plain
\newcommand{\Eachtplain}
{
\setlength{\unitlength}{0.6mm}
\begin{picture}(70,20)(-5,-15)

\put (0,0) {\circle{3}}

\put (2,0) {\line(1,0){6}}

\put (10,0) {\circle{3}}

\put (12,0) {\line(1,0){6}}

\put (20,0) {\circle{3}}

\put (22,0) {\line(1,0){6}}

\put (30,0) {\circle{3}}

\put (32,0) {\line(1,0){6}}

\put (40,0) {\circle{3}}

\put (42,0) {\line(1,0){6}}

\put (50,0) {\circle{3}}

\put (52,0) {\line(1,0){6}}

\put (60,0) {\circle{3}}

\put (20,-8) {\line(0,1){6}}

\put (20,-10) {\circle{3}}

\end{picture}
}

% Dn - plain
\newcommand{\Dnplain}
{
\setlength{\unitlength}{0.6mm}
\begin{picture}(70,20)(-5,-15)

\put (0,0) {\circle{3}}

\put (2,0) {\line(1,0){6}}

\put (10,0) {\circle{3}}

\put (12,0) {\line(1,0){6}}

\put (20,0) {\circle{3}}

\put (22,0) {\line(1,0){6}}

\put (32,0) {\circle*{1}}
\put (35,0) {\circle*{1}}
\put (38,0) {\circle*{1}}

\put (42,0) {\line(1,0){6}}

\put (50,0) {\circle{3}}

\put (52,0) {\line(1,0){6}}

\put (60,0) {\circle{3}}

\put (10,-8) {\line(0,1){6}}

\put (10,-10) {\circle{3}}

\end{picture}
}

% An - plain
\newcommand{\Anplain}
{
\setlength{\unitlength}{0.6mm}
\begin{picture}(70,10)(-5,-5)

\put (0,0) {\circle{3}}

\put (2,0) {\line(1,0){6}}

\put (10,0) {\circle{3}}

\put (12,0) {\line(1,0){6}}

\put (20,0) {\circle{3}}

\put (22,0) {\line(1,0){6}}

\put (32,0) {\circle*{1}}
\put (35,0) {\circle*{1}}
\put (38,0) {\circle*{1}}

\put (42,0) {\line(1,0){6}}

\put (50,0) {\circle{3}}

\put (52,0) {\line(1,0){6}}

\put (60,0) {\circle{3}}

\end{picture}
}

% roots
\newcommand{\roots}
{
\setlength{\unitlength}{1.5mm}
\begin{picture}(30,27)(-15,0)

\put (0,4) {\circle{3}}

\put (-1,0) {\makebox{$\alpha_{3,2}$}}

\put (2,6) {\line(2,1){4}}

\put (8,9) {\circle{3}}

\put (10,8) {\makebox{$\alpha_{1,2}$}}

\put (8,11) {\line(0,1){4}}

\put (8,17) {\circle{3}}

\put (10,16) {\makebox{$\alpha_{1,3}$}}

\put (-6,7) {\line(2,-1){4}}

\put (-8,9) {\circle{3}}

\put (-14,8) {\makebox{$\alpha_{3,1}$}}

\put (-8,11) {\line(0,1){4}}

\put (-8,17) {\circle{3}}

\put (-14,16) {\makebox{$\alpha_{2,1}$}}

\put (3,12) {\line(2,-1){3}}

\put (-6,16) {\line(2,-1){3}}

\put (2,21) {\line(2,-1){4}}

\put (-6,19) {\line(2,1){4}}

\put (3,14) {\line(2,1){3}}

\put (-6,10) {\line(2,1){3}}

\put (0,15) {\line(0,1){5}}

\put (0,6) {\line(0,1){5}}

\put (0,22) {\circle{3}}

\put (-1,25) {\makebox{$\alpha_{2,3}$}}

\put (0,13) {\circle{3}}

\end{picture}
}

% toric description of blow-up
\newcommand{\toric}{
\setlength{\unitlength}{.3mm}
\begin{picture}(110,110)(-10,-10)

\put (0,0) {\vector(1,0){90}}

\put (0,0) {\vector(0,1){90}}

\put (0,0) {\line(1,2){40}}

\multiput (0,0) (15,15) {4} {\line(1,1){12}}

\put (30,0) {\circle*{2}}

\put (30,30) {\circle*{2}}

\put (30,60) {\circle*{2}}

\put (33,-10) {\makebox{$x \equiv (1,0)$}}

\put (33,21) {\makebox{$y \equiv (1,1)$}}

\put (33,55) {\makebox{$z \equiv (1,2)$}}

\end{picture}
}

% cone
\newcommand{\cone}
{
\setlength{\unitlength}{0.5mm}
\begin{picture}(75,80)(-35,-65)

\cbezier(-20,0)(-20,5.522847498)(-11.045694996,10)(0,10)
\cbezier(0,10)(11.045694996,10)(20,5.522847498)(20,0)
\cbezier(20,0)(20,-5.522847498)(11.045694996,-10)(0,-10)
\cbezier(0,-10)(-11.045694996,-10)(-20,-5.522847498)(-20,0)

\put (-30,-53) {\vector(1,1){60}}
\put (32,5) {\makebox{$x$}}

\put (30,-53) {\vector(-1,1){60}}
\put (-35,5) {\makebox{$z$}}

\put (0,-23) {\vector(1,0){35}}
\put (30,-30) {\makebox{$y$}}

\cbezier[20](-20,-46)(-20,-40.477152502)(-11.045694996,-36)(0,-36)
\cbezier[20](0,-36)(11.045694996,-36)(20,-40.477152502)(20,-46)
\cbezier(20,-46)(20,-51.522847498)(11.045694996,-56)(0,-56)
\cbezier(0,-56)(-11.045694996,-56)(-20,-51.522847498)(-20,-46)

\end{picture}
}

% resolution - real image
\newcommand{\resreal}{
\setlength{\unitlength}{0.5mm}
\begin{picture}(170,105)(-150,-85)

% cone

\cbezier(-20,0)(-20,5.522847498)(-11.045694996,10)(0,10)
\cbezier(0,10)(11.045694996,10)(20,5.522847498)(20,0)
\cbezier(20,0)(20,-5.522847498)(11.045694996,-10)(0,-10)
\cbezier(0,-10)(-11.045694996,-10)(-20,-5.522847498)(-20,0)

\put (-30,-53) {\line(1,1){60}}

\put (30,-53) {\line(-1,1){60}}

\put (0,-23) {\circle{3}}

\cbezier[20](-20,-46)(-20,-40.477152502)(-11.045694996,-36)(0,-36)
\cbezier[20](0,-36)(11.045694996,-36)(20,-40.477152502)(20,-46)
\cbezier(20,-46)(20,-51.522847498)(11.045694996,-56)(0,-56)
\cbezier(0,-56)(-11.045694996,-56)(-20,-51.522847498)(-20,-46)

% cylinder

\cbezier[20](-100,-23)(-100,-17.477152502)(-91.045694996,-13)(-80,-13)
\cbezier[20](-80,-13)(-68.95430504,-13)(-60,-17.477152502)(-60,-23)
\cbezier(-60,-23)(-60,-28.522847498)(-68.95430504,-33)(-80,-33)
\cbezier(-80,-33)(-91.045694996,-33)(-100,-28.522847498)(-100,-23)

\put (-100,-53) {\line(0,1){60}}
\put (-60,-53) {\line(0,1){60}}

% morphism

\put (-51,-23) {\vector(1,0){34}}

% text

\put (-119,-65) {\makebox{
$\widetilde{X}_{\rm real} = T^{\ast} {\mathbb R}{\mathbb P}^1 $}}
\put (-100,-75) {\makebox{
$ = T^{\ast} S^1 = S^1 \times {\mathbb R}$}}

\put (-10,-65) {\makebox{$X_{\rm real}$}}

\put (-143,-23) {\makebox{$ E_{\rm real} = {\mathbb R}{\mathbb P}^1 $}}
\put (-127,-33) {\makebox{$=S^1$}}

\end{picture}
}

% ----------------------------------------------------------------------------------------------------------------------------------------------------------------------------------------------------------------------------------------
%
% settings and self-defined commands
%
% ----------------------------------------------------------------------------------------------------------------------------------------------------------------------------------------------------------------------------------------

\newcommand{\GG}{\ensuremath{\mathfrak g}}
\newcommand{\OO}{\mathcal O}
\newcommand{\Ox}{\ensuremath{\OO_X}}
\newcommand{\Oxx}{\ensuremath{ \OO_{\widetilde{X}}}}
\newcommand{\res}{\ensuremath{ \pi : \widetilde{X} \rightarrow X}}
\newcommand{\Znum}{\ensuremath{ Z_{\rm num}}}
\newcommand{\PP}{\ensuremath{{\mathbb P}^1}}
\newcommand{\CC}{\ensuremath{\mathbb C}}
\newcommand{\ZZ}{\ensuremath{\mathbb Z}}
\newcommand{\degz}{\text{deg}_{\Znum}}
\newcommand{\II}{\ensuremath{\mathcal J}}
\newcommand{\mm}{\mathfrak m}
\newcommand{\XX}{\ensuremath{\widetilde{X}}}
\newcommand{\XXX}{\ensuremath{\XX ^{\prime}}}

\pagestyle{plain}
\setlength{\parindent}{0mm}

\newtheorem{theorem}{Theorem}[section]
\newtheorem{proposition}[theorem]{Proposition}
\newtheorem{lemma}[theorem]{Lemma}
\newtheorem{corollary}[theorem]{Corollary}
\newtheorem{defprop}[theorem]{Definition and Proposition}
\newtheorem{remark}[theorem]{Remark}

\begin{titlepage}

\begin{flushright}
  Part III Essay \\ Easter 2002
\end{flushright}

\bigskip

\begin{center}
  \LARGE The Dynkin Diagrams of Rational Double Points
\end{center}

\bigskip

\begin{center}
  \normalsize Benjamin Friedrich \\
  \smallskip
  \sl Trinity College \\
  \bigskip
  {\it Home address:} \\
  Kleiner Warnowdamm 10 \\
  18109 Rostock \\
  Germany \\
  {\tt bf216@hermes.cam.ac.uk} \\
  \bigskip
  {\rm 15 May 2002}
\end{center}

\bigskip

\begin{center}
  I declare that this essay is work done as part of the Part III Examination. It is the
  result of my own work, and except where stated otherwise, includes nothing which was
  performed in collaboration. No part of this essay has been submitted for a degree or
  any such qualification. \\[2cm]
  Benjamin Friedrich
\end{center}

\end{titlepage}

\begin{titlepage}

\pagebreak
\thispagestyle{empty}

\begin{flushright}
  Part III Essay \\ Easter 2002
\end{flushright}

\bigskip

\begin{center}
  \LARGE The Dynkin Diagrams of Rational Double Points
\end{center}

\bigskip

\setlength{\parindent}{0mm}
\begin{abstract}

\begin{center}
Rational double points are the simplest surface singularities. In
this essay we will be mainly concerned with the geometry of the
exceptional set corresponding to the resolution of a rational
double point. We will derive the classification of rational
double points in terms of Dynkin diagrams.
\end{center}

\end{abstract}

\end{titlepage}

\goodbreak

\tableofcontents

\goodbreak

% ----------------------------------------------------------------------------------------------------------------------------------------------------------------------------------------------------------------------------------------
\section{Introduction}

Rational double points\footnote{In the literature, they are also
called Du Val singularities, Kleinian singularities or simple
critical points. } are the simplest surface singularities and were
first studied by Du Val \cite{DV}. One may think of them as
neglible singularities. They play an important r\^{o}le in the
classification of surfaces and occur in the theory of
simultaneous resolutions of singularities.

In this essay we will be mainly concerned with the geometry of the
exceptional set corresponding to a resolution of a rational
double point. We will derive the classification of rational
double points in terms of Dynkin diagrams. It should be noted,
that the proof of this classification is rather lengthy. However,
the author was unable to find the complete proof in a single
source and therefore decided to present it in full detail. Most
ideas are taken from two papers of Artin \cite{A2}, \cite{A1},
balanced with a slightly different approach in Reid's draft
\cite{Reid}. Further parts of the argument are taken from Durfee
\cite{Df}, Mumford \cite{Mf} and Brieskorn \cite{B3}. The second
article by Pinkham in \cite{DPT} treats the topic very nicely,
although some difficult steps are omitted.

Finally, we will find a connection between the most simple objects
in different fields of mathematics: Rational double points are
linked with Platonic solids and simple Lie groups.

% ----------------------------------------------------------------------------------------------------------------------------------------------------------------------------------------------------------------------------------------
\section{Basic facts on surface singularities}

% ----------------------------------------------------------------------------------------------------------------------------------------------------------------------------------------------------------------------------------------
\subsection{Definitions}

We want to study {\em surface singularities} $(X,x)$; here $X$ is a nor\-mal, two-dimen\-sional, projective variety over $\mathbb C$ which is non-singular,
except maybe at $x \in X$.

Two singularities are isomorphic, if there exist open neighbourhoods of the singular points which are isomorphic.

A {\em resolution} of $(X,x)$ is a birational, proper and
surjective morphism

$$ \res $$

where $\widetilde{X}$ is a non-singular projective variety over \CC.

See section \ref{blowup} for an example.

It is an important and difficult theorem, that resolutions always exist; for a general discussion we refer to \cite{Lp}, \cite{Lf}.

\paragraph{Immediate properties of the exceptional set}
The {\em exceptional set} $E:=\pi^{-1}(x)$ is compact (since $X$
is proper) and one-dimensional (since $\pi$ is birational).
Moreover it is connected by Zariski's connectedness theorem
\ref{app5}. Therefore $E$ is a bunch of irreducible projective
curves

$$E=\bigcup\limits_{i=1}^{n}E_i .$$

We say that a surface singularity is
\begin{description}
\item{\em rational}, if for a resolution
$$ \res $$
the first higher direct image sheaf of $\XX\text{'s}$ structure sheaf vanishes
\begin{equation}
\label{ratcond}
R^1 \pi_{\ast} \Oxx = 0,
\end{equation}
and
\item{\em a double point}, if the local ring $\OO _{X,x}$ has multiplicity two, i.e. the
leading coefficient of its Hilbert-Samuel polynomial is two
(\cite{Ng} III \S 23, \cite{ZS} vol. 2, VIII \S 10).
\end{description}

\paragraph{Remarks:}
\begin{enumerate}
\item
\begin{enumerate}
\item The definition of a rational singularity is independent of the chosen resolution:
Since $ R^1 \pi_{\ast} \Ox~ $ is a coherent sheaf (\cite{Ha} III.8.8.(b))
concentrated on $x$,
all we are interested in is $h^0(X,R^1 \pi_{\ast} \Oxx)$.
However, we will see soon in section \ref{leray} that
$$ p_a(X) - p_a(\widetilde{X}) = h^0(X,R^1 \pi_{\ast} \Oxx) $$
and the arithmetic genus of a {\bf non-singular} projective variety is a birational invariant (\cite{Ha} V.5.6).

\item The condition (\ref{ratcond}) may appear opaque at a first glance, but will hopefully become more transparent in the sequel.
For example, it implies that the $E_i$ are rational curves.

\end{enumerate}

\item Since we are in the normal case, the condition for a double point means, that two general curves on $X$ through $x$ have local intersection number
two at $x$ (\cite{Reid} 4.6). If $X$ is a hypersurface $f^{-1}(0)$, yet another way to state this condition is

$$ f \in m_x^2 \text{ and } f \notin m_x^3 ,$$

where $m_x$ is the ideal of functions vanishing at $x$ (\cite{Dm} (7.48)).
\end{enumerate}

% ----------------------------------------------------------------------------------------------------------------------------------------------------------------------------------------------------------------------------------------------
\subsection{A first consequence of the rationality condition (\ref{ratcond})}
Let us mention a simple consequence of the rationality condition (\ref{ratcond}).

\begin{proposition}
Let $\res$ be a resolution as above. If
$$ R^1 \pi_{\ast} \Oxx = 0 $$
then
$$ p_a(X) = p_a(\widetilde{X}) . $$
\end{proposition}

We need a lemma.

\begin{lemma}
For any resolution $\res$, we have
$$ \pi_{\ast} \Oxx = \Ox . $$
\end{lemma}

{\bf Proof} (\cite{Ha} p. 280):
Since the question is local on $X$, we can assume $X$ is affine, say $X=\text{Spec } A$.
By \cite{Ha} II.5.8.(b), $\pi_{\ast} \Oxx$ is a coherent sheaf of $\Ox$-algebras, hence $B:=H^0(X,\pi_{\ast} \Oxx)$ is a finitely generated $A$-module.
But $A$ and $B$ are integral domains with the same quotient field (since $\pi$ is birational)
and $A$ is algebraically closed (since $X$ is normal), thus $B=A$, and $ \pi_{\ast} \Oxx = \Ox$. \hfill $\Box$

\vspace{3mm} {\bf Proof of the proposition} (\cite{Ha} Ex. III.8.1):
Let
\begin{equation}
\label{injres}
0 \longrightarrow I^0 \stackrel{d^0}{\longrightarrow} I^1 \stackrel{d^1}{\longrightarrow} I^2 \longrightarrow \dots
\end{equation}
be an injective resolution for $\Oxx$. We have not only $R^1
\pi_{\ast} \Oxx = 0$, but $R^i \pi_{\ast} \Oxx = 0$ for $i \ge 1$,
because the fibers of $\pi$ have dimension $\le 1$  (\ref{app7}).
Therefore, by applying $\pi_{\ast}$ to (\ref{injres}), we obtain
again an exact sequence
\begin{equation}
\label{acres}
0 \longrightarrow \pi_{\ast}\Oxx \longrightarrow \pi_{\ast}I^0 \stackrel{\pi_{\ast}d^0}{\longrightarrow} \pi_{\ast}I^1 \stackrel{\pi_{\ast}d^1}
{\longrightarrow} \pi_{\ast}I^2 \longrightarrow \dots .
\end{equation}
Since injectives are flasque (\cite{Ha} III.2.4), direct images of flasque sheaves are flasque, and flasque sheaves are acyclic for the global section
functor (\cite{Ha} III.2.5), we see that (\ref{acres}) is an acyclic resolution for $\Ox = \pi_{\ast} \Oxx$.
Thus
\begin{eqnarray*}
H^i(X,\Ox) & = & \frac{\text{ker }H^0(X,\pi_{\ast}d^i)}{\text{im }H^0(X,\pi_{\ast}d^{i-1})} \\
& = & \frac{\text{ker }H^0(\widetilde{X},d^i)}{\text{im }H^0(\widetilde{X},d^{i-1})} \\
& = & H^i(\widetilde{X},\Oxx)
\end{eqnarray*}
and our claim follows. \hfill $\Box$

\paragraph{Remark:}
\label{leray}
This proof is just a degenerated case of the {\em Leray spectral sequence}
\begin{eqnarray*}
&& E_2^{p,q} = H^p(X,R^q\pi_{\ast}\Oxx) \\
&& \Rightarrow \\
&& E_{\infty}^{p,q}=H^{p+q}(\widetilde{X},\Oxx)
\end{eqnarray*}
(cf. \cite{Gd} II.4.17.1, \cite{Wb} V; for general information on spectral sequences cf. \cite{McC}, \cite{Sh} II \S  4),
which takes in our setting the simple form (theorem \ref{app6} and \ref{app7})

\goodbreak

$$ E_2^{p,q}:
\begin{array}{cccc}
0 & 0 & 0 & 0 \\
H^0(X,R^1\pi_{\ast}\Oxx) & ? & ? & 0 \\
H^0(X,\Ox) & H^1(X,\Ox)  & H^2(X,\Ox) & 0
\end{array}
$$
\nopagebreak[4] \hrulefill \nopagebreak[4]
\begin{eqnarray*}
& H^i(\widetilde{X},\Oxx)=E_{\infty}^{i,0}=E_2^{i,0}=H^i(X,\Ox) \text{ for } i=0,1, \\
& H^2(\widetilde{X},\Oxx)=E_{\infty}^{2,0}=E_2^{2,0}/E_2^{0,1}=H^2(X,\Ox)/H^0(X,R^1\pi_{\ast}\Oxx).
\end{eqnarray*}
From this we get
$$ p_a(X) - p_a(\widetilde{X}) = h^0(X,R^1 \pi_{\ast} \Oxx) $$
and our proposition (and its converse!) follow at once.

% -----------------------------------------------------------------------------------------------------------------------------------------------------------------------------------------------------------------------------------------------
\subsection{Further properties of the exceptional set $E$}

It will be a great technical convenience to work with {\em good resolutions}; for them we require that

\begin{enumerate}

\item all $E_i$ are non-singular,

\item \label{goodcond2}
$E_i \cap E_j \cap E_k = \varnothing $ for mutually distinct $i, j, k$,

\item the intersection of $E_i$ and $E_j$ is transverse for $i \ne j$.

\end{enumerate}

Any resolution $\res$ of a {\bf surface} $X$ can be brought in
such a nice form by successively blowing up points of
$\widetilde{X}$ (cf. again \cite{Lp}, \cite{Lf}, and also
\cite{Ha} V.3.8, V.3.9).

In the following we do always assume that
$ \res $
is good.

A fundamental fact about good resolutions is the following

\begin{proposition} {\rm (\cite{Mf} p. 6)}
\label{negdef}
The {\em intersection matrix of the resolution} $\left( E_i \cdot E_j \right)_{i,j=1 \dots n}$ is negative definite.
\end{proposition}

{\bf Proof}:
We take a meromorphic function $f \in k(X)$ with $f(x)=0$ and define two effective divisors

$$H_0:=f^{-1}(0) \text{ and } H_{\infty}:=f^{-1}(\infty) . $$

Denote the proper transform of $H_i$ with $\widetilde{H}_i$ for $i=0, \infty$ respectively.
Then we have a linear equivalence of divisors

$$ \widetilde{H}_{\infty} \sim \widetilde{H}_0 + \sum\limits_{i=1}^n m_i E_i$$

where $m_i=\text{ord} _{E_i} f \circ \pi > 0$.

It suffices to prove that the matrix

$$M:=\left(m_i E_i \cdot m_j E_j\right)_{i,j=1, \dots ,n}$$

is negative definite. Now we have
$M_{i,j} \ge 0$ if $i \ne j$ (since the $E_i$ are irreducible)
and
\begin{eqnarray*}
\sum\limits_{i=1}^n M_{i,j} & = & \sum\limits_{i=1}^n (m_i E_i \cdot m_j E_j) \\
& = & (\widetilde{H}_{\infty} - \widetilde{H}_0) \cdot m_j E_j \\
& = & 0 - \widetilde{H}_0 \cdot m_j E_j \le 0 .
\end{eqnarray*}

This implies that $M$ is negative semi-definite:

\begin{eqnarray}
\sum\limits_{i,j=1}^n a_i a_j M_{i,j} & = & \sum\limits_{i=1}^n a_i^2 M_{i,i} + 2\sum\limits_{i,j=1 \atop i<j}^n a_i a_j M_{i,j} \nonumber\\
& = & \sum\limits_{j=1}^n
\underbrace{
\left(
\sum\limits_{i=1}^n M_{i,j}
\right)
}_{\le 0}
a_j^2
\label{mf1}
- \sum\limits_{i,j=1 \atop i<j}^n \underbrace{M_{i,j} (a_i-a_j)^2}_{\ge 0} \le 0 .
\end{eqnarray}

To show definiteness, we note that $\widetilde{H}_0$ must pass through some $E_i$, hence
$$ \sum \limits_{i=1}^n M_{i,j_0} < 0 \text{ for some } j_0 . $$
Suppose we have equality in (\ref{mf1}). Then $a_{j_0}=0$.
Furthermore, we get $a_i=a_j$ if $M_{i,j} > 0$, or inductively $a_i=a_j$ if $E_i$ and $E_j$
are connected in $E$. But $E$ is connected, hence $a_i=0$ for $i=1, \dots , n$ in this case. \hfill $\Box$

\vspace{3mm}
In the proof of proposition \ref{negdef}, we encountered an effective {\em exceptional divisor} (i.e. an divisor supported on $E$)
$$ Z:=\widetilde{H}_{\infty}-\widetilde{H}_0=\sum \limits_{i=1}^n m_i E_i > 0 $$
which had the note-worthy property
\begin{equation}
\label{numcond}
Z \cdot E_i \le 0 \text{ for } i=1, \dots , n.
\end{equation}

Since $E$ is connected, any exceptional divisor $Z$ with this property (\ref{numcond}) must satisfy $Z \ge E$, by the arguments used in that proof. If two
exceptional divisors
$Z^1=\sum \limits_{i=1}^n r_i^1 E_i > 0$ and $Z^2=\sum \limits_{i=1}^n r_i^2 E_i > 0$ both satisfy (\ref{numcond}) then obviously so does
$Z:=\min(Z^1,Z^2)=\sum \limits_{i=1}^n r_i E_i > 0$ where $r_i:=\min(r_i^1,r_i^2)$:
$
Z \cdot E_i \le Z^j \cdot E_i \le 0$ whenever $r_i=r_i^j
$.
Hence there exists a minimal positive exceptional divisor, called
the {\em numerical divisor} $\Znum$ (\cite{Reid} 4.5) (also
called {\em fundamental divisor} \cite{A1}), for which
(\ref{numcond}) holds.

This divisor $\Znum$ will provide a useful tool to describe the exceptional set of a rational singularity.

% ------------------------------------------------------------------------------------------------------------------------------------------------------------------------------------------------------------------------------------------
\section{The geometry of the exceptional set $E$ of a resolution of a rational singularity}

Throughout this section, we assume that \res\ is good resolution
of a rational singularity $(X,x)$ and $E$ its exceptional set.

We will prove that the $E_i$ are rational curves $E_i \cong {\mathbb P}^1$.
Moreover, we will be able to read off the multiplicity of $(X,x)$ as the self-intersection-number $-(\Znum)^2$.

The idea is to study fatter and fatter infinitesimal neighbourhoods of $E$ in order to examine the embedding of $E$ in $\widetilde{X}$.

We will identify an exceptional divisor $Z=\sum\limits_{i=1}^n r_i E_i$ with its associated {\em positive cycle}: this is the, generally {\bf non-reduced},
scheme $(\text{Supp } Z, \OO _Z)$.
Recall that $\OO_Z = \text{coker} (\Oxx (-Z) \rightarrow \Oxx )$, i.e. $(\text{Supp } Z, \OO _Z)$ is the subscheme of $\widetilde{X}$ defined by the coherent
sheaf of ideals on $\widetilde{X}$ whose sections on an open $U \subset \widetilde{X}$ are the
 rational functions $f \in \Gamma(U,\Oxx)$ which have zeros of order at
least $r_i$ along $E_i$ for all $i$ with $E_i \cap U \neq \varnothing $ \cite{A2}.
Note $\text{Supp } Z = \bigcup\limits_{r_i > 0} E_i$.

% ----------------------------------------------------------------------------------------------------------------------------------------------------------------------------------------------------------------------------------------
\subsection{The exceptional curves $E_i$ are rational}

\begin{theorem}
\label{theo1}
{\rm (\cite{A1} prop. 1, \cite{B3} Lemma 1.3)}
The exceptional set of a good resolution of a rational singularity
consists of rational projective curves $E_i \cong {\mathbb P}^1$.
\end{theorem}

{\bf Proof} (\cite{A1}, \cite{A2}):

The proof relies on Grothendieck's theorem on formal functions \ref{app3}, which takes in our case the form

$$ 0 = \left( \left( R^1 \pi_{\ast} \Oxx \right)_x \right)\widehat{}\ \cong\
  \mathop{\varprojlim}_{k=1}^{\infty} H^1\left( E, \Oxx \otimes _{\Ox} {\OO_{X,x} / {m_x}^k} \right) $$

where $m_x \subset \OO_{X,x}$ is the maximal ideal corresponding to $x$ and completion is taken with respect to the $m_x$-adic topology.

(We will see later in lemma \ref{mznum}, that $\Oxx \otimes _{\Ox} {\OO_{X,x} / {m_x}^k } = \OO_{k \Znum}$.)

Since $E$ is one-dimensional, the natural map
$$ \Oxx \otimes _{\Ox} {\OO_{X,x} / {m_x}^{k+1} } \twoheadrightarrow \Oxx \otimes _{\Ox} {\OO_{X,x} / {m_x}^k }$$
induces a surjection on cohomology
$$ H^1\left( E, \Oxx \otimes _{\Ox} {\OO_{X,x} / {m_x}^{k+1} } \right) \twoheadrightarrow
  H^1\left( E, \Oxx \otimes _{\Ox} {\OO_{X,x} / {m_x}^k } \right)$$
by a vanishing theorem of Grothendieck \ref{app6}.
Thus we see
$$ H^1 \left( E, \Oxx \otimes _{\Ox} {\OO_{X,x} / {m_x}^{k+1} } \right) = 0 \text{ for all } k \in {\mathbb N}. $$

We denote the sheaf of ideals of functions vanishing at $x$ by $\mm _x$.
Clearly every function in $\mm _x \cdot \Oxx$ vanishes on $E$; hence for every positive cycle $Z$ we can find an integer $k$ such that every
function in ${\mm _x}^k \cdot \Oxx$ vanishes on $Z$. We now have
$$ \Oxx \otimes _{\Ox} {\OO_{X,x} / {m_x}^k } \twoheadrightarrow \OO _Z $$
and (\ref{app6})
$$ 0= H^1 \left( E, \Oxx \otimes _{\Ox} {\OO_{X,x} / {m_x}^{k+1} } \right) \twoheadrightarrow H^1(E, \OO _Z). $$

In particular $H^1(E, \OO _{E_i}) =0$ for $i=1,\dots,n$, from which we conclude
$p_a(E_i) \nolinebreak = \nolinebreak 0$, i.e. $E_i \cong {\mathbb P}^1$. \hfill $\Box$

\begin{corollary}
\label{H10}
In the proof of theorem \ref{theo1} we have just seen $H^1(E, \OO _Z)=0$ for every positive cycle $Z$.
\end{corollary}

We can make a more precise statement for the numerical divisor \Znum.

\begin{corollary}
\label{znum0}
{\rm (\cite{A1} thm. 3)} With the assumptions of the
theorem we have $$ p_a(\Znum) =0 . $$
\end{corollary}

{\bf Proof} (\cite{Reid} 3.11):
The statement follows immediately from corollary \ref{znum0}
and the general fact
\fbox{$h^0(E,\OO_{\Znum}) = 1$}.
This can be proved easily by induction:
We know $h^0(E)=1$. Assume $h^0(Y)=1$ for a positive cycle
$E \le Y \lneq \Znum$.
We have
$$ Y \cdot E_i > 0 \text{ for some } i, $$
or equivalently $\text{deg}_{E_i} \OO_E(-Y) \le -1$,
by the very definition of \Znum. Certainly $Y+E_i \le \Znum$ in this situation.
From
$$
0=H^0(E_i,\OO_{E_i}(-Y)) \rightarrow H^0(E,\OO_{Y+E_i}) \rightarrow H^0(E,\OO_Y)
$$
we conclude
$h^0(E,\OO_{Y+E_i})=1$. \hfill $\Box$

% ----------------------------------------------------------------------------------------------------------------------------------------------------------------------------------------------------------------------------------------
\subsection{A criterion for rationality}

\begin{theorem}
\label{theo2}
{\rm (\cite{A1}, thm. 3)}

Conversely, if we have

$$p_a(\Znum) = 0$$

for the numerical cycle of a good resolution of a singularity $(X,x)$, then $(X,x)$ is rational.
\end{theorem}

We need the following lemma.

\begin{lemma}
Let $Z=\sum\limits_{i=1}^n r_i E_i$ be a positive cycle with the property that
$p_a(Y) \le 0$ for all positive cycles $Y \le Z$. Then $H^1(E, \OO _Z)=0$.
\end{lemma}

{\bf Proof}:
In particular $p_a(E_i)=0$ for all $i$ with $r_i \ge 1$, i.e. $E_i \cong {\mathbb P}^1$.
We use induction on $\sum\limits_{i=1}^n r_i$. Assume $H^1(E,\OO_Z) \neq 0$.
Let $Z_i:=Z-E_i$ for $r_i \ge 1$. By induction hypothesis
$$H^1(E, \OO _{Z_i})=0 , $$
hence for the kernel $M$
$$ 0 \rightarrow M \rightarrow \OO _Z \rightarrow \OO _{Z_i} \rightarrow 0 $$
we get (\ref{app6})
$$ H^1(E,M) \rightarrow H^1(E, \OO _Z) \rightarrow 0 , $$
i.e. $H^1(E,M) \neq 0$.
By the snake-lemma
$$
\begin{CD}
{}   @.   {}   @.   0   @>>>   M   @.   {} \\
@.        @.         @VVV        @VVV    @. \\
0 @>>> \Oxx(-Z) @>>> \Oxx @>>> \OO _Z @>>> 0 \\
@.        @VVV    @VVV        @VVV    @. \\
0 @>>> \Oxx(-Z_i) @>>> \Oxx @>>> \OO _{Z_i} @>>> 0 \\
@.        @VVV    @VVV       @.          @. \\
{} @.  \Oxx(-Z_i) \otimes \OO _{E_i} @>>> 0 @. {} @. {}
\end{CD}
$$

we deduce $M  \cong \Oxx(-Z_i) \otimes \OO _{E_i}$. Hence we can write
$$ 0 \neq H^1(E, \Oxx (-Z_i) \otimes \OO _{E_i}) = H^1(E_i,\Oxx(-Z_i) \otimes \OO _{E_i}) . $$
But $E_i$ is just ${\mathbb P}^1$, thus (\cite{Ha} III.5.1)
$$ \text{deg } \Oxx(-Z_i) \otimes \OO _{E_i} \le -2 . $$
On the other hand
$$ \text{deg } \Oxx(-Z_i) \otimes \OO _{E_i} = - Z_i \cdot E_i$$
and by the adjunction formula \ref{app2} we get
$$ Z \cdot E_i = (Z_i + E_i) \cdot E_i \ge 2 + {E_i}^2 = -K \cdot E_i. $$
Summing up $(Z+K) \cdot E_i \ge 0$ yields with the adjunction formula \ref{app2}
$$ 2 p_a(Z) - 2 = ( Z+K) \cdot Z \ge 0, $$
i.e. $p_a(Z) \ge 1$, a contradiction. \hfill $\Box$

\vspace{3mm}

{\bf Proof of the theorem} (\cite{A1} prop. 1, thm. 3):
We have seen in the proof of theorem \ref{theo1}, that by Grothendieck's theorem on formal functions (\ref{app3})
$$ 0 = \left( \left( R^1 \pi_{\ast} \Oxx \right)_x \right)\widehat{}\ \cong\
\mathop{\varprojlim}_{k=1}^{\infty} H^1( E, \OO _{k\Znum}), $$
hence it is sufficient to prove
$$ H^1(E, \OO _{k \Znum}) = 0 \text{ for all } k. $$
We already know $H^1(E, \OO _{\Znum})=0$ (since $p_a(\Znum )=0$ and $h^0(E,\OO_{\Znum})=1$).
From the surjection (\ref{app6})
$$
H^1(E, \OO _{\Znum}) \twoheadrightarrow H^1(E, \OO _{E_i})
$$
we find that $p_a(E_i)=0$, i.e. $E_i \cong {\mathbb P}^1$.

By the lemma, it is enough to show $p_a(Y) \le 0$ for all positive cycles $Y$. Let $Y_1:=Y$ and define $Y_{n+1}$ inductively as follows

\begin{enumerate}
\item if $Y_n \ge \Znum$, then $Y_{n+1}:=Y_n-\Znum \ge 0$.
\item if $Y_n \not\ge \Znum$, then $Y_n \cdot E_i > 0$ for some $i$ by the definition of \Znum. Choose such an $i$ with smallest possible multiplicity in
$Y_n$ and set $Y_{n+1}:=Y_n + E_i$.
\end{enumerate}
Stop when $Y_n=0$. We use the equation (\cite{Ha} Ex. V.1.3)
\begin{equation}
\label{lichtenbaum}
p_a(Z_1 + Z_2) = p_a(Z_1) + p_a(Z_2) + Z_1 \cdot Z_1 -1
\end{equation}
to calculate the arithmetic genus:

In case 1:
\begin{eqnarray*}
p_a(Y_n) & = & p_a(Y_{n+1}+\Znum) \\
& = & p_a(Y_{n+1}) + p_a(\Znum) + Y_{n+1} \cdot \Znum -1 \\
& \le & p_a(Y_{n+1})-1
\end{eqnarray*}

In case 2:
\begin{eqnarray*}
p_a(Y_{n+1}) & = & p_a(Y_n) + p_a(E_i) + Y_n \cdot E_i -1 \\
& \ge & p_a(Y_n).
\end{eqnarray*}

Steps of type $(2)$ cannot be repeated infinitely often without reaching a stage where $Y_n \ge \Znum$. Using equation (\ref{lichtenbaum}) once again, we see
that $p_a(Y)$ is a quadratic form in the coefficients $s_i$ of $Y=\sum\limits_{i=1}^n s_i E_i$ whose quadratic term is
$\frac{1}{2} \sum\limits_{i,j=1}^n s_i s_j E_i \cdot E_j$. But the matrix $(E_i \cdot E_j)_{i,j=1,\dots,n}$
is negative definite, hence $p_a(Y)$ is bounded above.
Consequently, there
can be only a finite number of steps and the algorithm must terminate with $Y_n=0$. Then $Y_{n-1}=\Znum$ and we have
$$ p_a(Y_1) \le \dots \le p_a(Y_{n-1})=0 . $$
\hfill $ \Box $

% ----------------------------------------------------------------------------------------------------------------------------------------------------------------------------------------------------------------------------------------
\subsection{Invertible sheaves on a positive cycle $Z$}

It is a natural question to ask what are the invertible sheaves on
a positive cycle $Z$, i.e. what is $\text{Pic } Z$?

In our particular case the answer is quite simple. It will provide an important tool in exploring the geometry of $E$ further.

For any invertible sheaf $\mathcal F$ on the positive cycle $Z \ge E$, we can define its {\em multidegree}

$$ \text{deg}_Z: \text{Pic }Z \rightarrow {\mathbb Z}^n $$

via the composite maps
$$
\text{Pic } Z \rightarrow \text{Pic } E_i \stackrel{\text{deg}}{\longrightarrow} {\mathbb Z}
\text{ for } i=1, \dots , n .
$$

Using {\em local transverse cuts} it is easily seen that this map is surjective
$ \text{deg}_Z \nolinebreak \text{Pic } \nolinebreak Z \nolinebreak \twoheadrightarrow \nolinebreak {\mathbb Z}^n $:
Choose a general point $p$ on any $E_i$ and construct a Cartier divisor
$\{ (U_j, f_j ) \}$ with support
$p$
and degree $1$ on $E_i$ whose local equation $s \in \OO _{Z,p}$ restricts to a local equation of $p$ in $\OO _{E_i, p}$. This gives
$$ \text{deg} _Z \{ (U_j, f_j ) \} = (0, \dots ,0 , \mathop{1}_i, 0, \dots ,0 ) . $$

In fact, we will prove that $ \text{deg}_Z \text{Pic } Z \twoheadrightarrow {\mathbb Z}^n $ is an isomorphism.

It is a well-known fact that (\cite{Ha} Ex. III.4.5) $$ \text{Pic
} Z = H^1(E, \OO _Z ^{\ast} ). $$ (One may think of a
\v{C}hech-$1$-cocycle $\{ ( U_i \cap U_j , g_{i,j} ) \} \in H^1(E,
\OO_Z^{\ast})$ as a set of transition functions $g_{i,j} :
\OO_{U_i \cap U_j} \stackrel{\sim}{\longrightarrow} \OO_{U_i \cap
U_j} $ which define a line bundle on $Z$.)

In the reduced case $Z=E$, it is easy to see what $H^1(E, \OO_E^{\ast})$ is, if
we allow transcendental methods. Let $\circ_h$ denote the functor
form the category of schemes of finite type over \CC\ to the category of complex analytic spaces.
(cf. \cite{Ha} B and section \ref{complexana}).

Since $E$ is projective over \CC, a theorem by Serre (\cite{Ha}
B.2.1) tells us that
$$
H^i(E,{\mathcal F}) \cong H^i(E_h,{\mathcal F}_h)
$$
for every coherent sheaf $\mathcal F$ on $E$.
The exponential sequence (\cite{Ha} V.5)
$$
0 \rightarrow \ZZ \rightarrow \OO_{E_h} \stackrel{\text{exp} 2 \pi i }{\longrightarrow} \OO_{E_h}^{\ast} \rightarrow 0
$$
yields (corollary \ref{H10} and theorem \ref{app6})

$$ 0 \rightarrow H^1(E_h, \OO_{E_h}^{\ast}) \rightarrow H^2(E_h, \ZZ) \rightarrow 0, $$

i.e.

$$ \text{Pic } E \cong H^2(E_h, \ZZ). $$

As we have already seen, E is built up out of $n$ spheres $S^2
\cong \PP$, which intersect each other transversely. Hence by the
Mayer-Vietoris sequence for, say singular cohomology

$$ H^2(E_h, \ZZ) \cong \ZZ^n . $$

(We will see later, that $E$ has the homotopy type of a bouquet of $n$ spheres $E \simeq (S^2)^{\vee n}$.)
Therefore, we obtain
\begin{equation*}
\text{Pic } E \cong \ZZ^n .
\end{equation*}

Unfortunately, there is no analogue of the exponential sequence in
the non-reduced case. Instead, we need the following proposition
by Artin, whose proof uses a "first order exponential".

\begin{proposition}
\label{picz}
{\rm (\cite{A2} lemma 1.4)}
We have
$$
H^1(E,\OO_Z) \cong H^1(E,\OO_E) \cong {\mathbb Z}^n
$$
for every positive cycle $Z \ge E$.
\end{proposition}

{\bf Proof}: We will proceed by induction: The case $Z=E$ is
trivial, so assume the proposition holds for
$Z^{\prime} =  Z - E_i \ge E$. We fix our notation for the
following kernels

\newcommand{\Nc}{\mathcal N}
\newcommand{\Mc}{\mathcal M}
\newcommand{\Ic}{\mathcal J}
\newcommand{\Kc}{\mathcal K}

\begin{eqnarray*}
& 0 \rightarrow \Nc \rightarrow \OO_Z \rightarrow
\OO_E \rightarrow 0 , \\
& 0 \rightarrow \Mc \rightarrow \OO_Z^{\ast} \rightarrow
\OO_E^{\ast} \rightarrow 0 , \\
& 0 \rightarrow \Nc ^{\prime} \rightarrow \OO_{Z^{\prime}} \rightarrow
\OO_E \rightarrow 0 , \\
& 0 \rightarrow \Mc ^{\prime} \rightarrow \OO_{Z^{\prime}}^{\ast} \rightarrow
\OO_E^{\ast} \rightarrow 0 , \\
& 0 \rightarrow \Ic \rightarrow \OO_Z \rightarrow
\OO_{Z^{\prime}} \rightarrow 0 & \text{ and } \\
& 0 \rightarrow \Kc \rightarrow \OO_Z^{\ast} \rightarrow
\OO_{Z^{\prime}}^{\ast} \rightarrow 0 .
\end{eqnarray*}
By \ref{app6}, it suffices to prove $H^1(E,\Mc)=0$.

Note that $H^0(E,\OO_E) = \CC$ (and also
$H^0(E,\OO_E^{\ast})=\CC^{\ast}$), since $E$ is connected. In
particular, we get a surjection
$$
H^0(E,\OO_Z) \twoheadrightarrow H^0(E,\OO_E) ,
$$
which implies (corollary \ref{H10})
$$
H^1(E,\Nc)=0 .
$$
Similarly, $H^1(E,\Mc ^{\prime})=0$ (using the induction
hypothesis). Now, these kernels are linked by the short exact
sequences
\begin{eqnarray*}
& 0 \rightarrow \Ic \rightarrow \Nc \rightarrow \Nc ^{\prime}
\rightarrow 0 \\
& 0 \rightarrow \Kc \rightarrow \Mc \rightarrow \Mc ^{\prime}
\rightarrow 0
\end{eqnarray*}
and we obtain
\begin{eqnarray*}
& H^0(E,\Nc ^{\prime}) \stackrel{\delta}{\longrightarrow}
H^1(E,\Ic) \rightarrow 0 \\
& H^0(E,\Mc ^{\prime})
\stackrel{\delta ^{\prime}}{\longrightarrow}
H^1(E,\Kc) \rightarrow H^1(E, \Mc) \rightarrow 0 .
\end{eqnarray*}

Because of $\Ic \cdot \Nc = 0$ (thus $\Ic ^2 = 0$), we have an
isomorphism
$$
\epsilon : \Ic \stackrel{\sim}{\longrightarrow} \Kc
$$
via
$$
s \in \Gamma(U,\Ic) \mapsto 1+s \in \Gamma(U,\Kc),
$$
hence $H^1(E,\Ic) \cong H^1(E, \Kc)$.
Analogously, we have a bijection (not a morphism, in general!)
\begin{eqnarray*}
\epsilon ^{\prime} : H^0(E, \Nc ^{\prime})
& \rightarrow &
H^0(E, \Mc ^{\prime}) \\
s^{\prime} & \mapsto & 1+s^{\prime} .
\end{eqnarray*}
Therefore, it suffices to show that the following diagram commutes
$$
\begin{CD}
H^0(E,\Nc ^{\prime}) @>{\delta}>> H^1(E, \Ic) \\
@V{\epsilon^{\prime}}VV @VV{\epsilon}V \\
H^0(E,\Mc ^{\prime}) @>{\delta^{\prime}}>> H^1(E, \Kc) .
\end{CD}
$$

Pick an element $s^{\prime} \in H^0(E,\Nc ^{\prime})$ and choose
an open covering $\{U_i\}$ of $E$ such that $s^{\prime}$ can be
lifted to $s_i \in \Gamma(U_i,\Nc)$. Now we can write
$\delta(s^{\prime})$ as the \v{C}ech-$1$-cocycle
$$
\{(U_i \cap U_j, s_i - s_j )\} \in H^1(E, \Ic)
$$
and get
$$
\epsilon(\delta(s^{\prime})) =
\{(U_i \cap U_j, 1+ s_i - s_j )\} \in H^1(E, \Kc) .
$$
In the same way, we can lift
$\epsilon^{\prime}(s^{\prime})=1+s^{\prime}$
to $1+s_i \in \Gamma(U_i, \Mc)$ and obtain
$$
\epsilon^{\prime}(\delta^{\prime}(s^{\prime})) =
\{(U_i \cap U_j, \frac{1 + s_i}{1 + s_j} )\} \in H^1(E, \Kc) .
$$
We use $\Ic \cdot \Nc$ in order to show
$$
\epsilon(\delta(s^{\prime})) =
\epsilon^{\prime}(\delta^{\prime}(s^{\prime})) .
$$
Since
$$
s_i - s_j \in \Gamma(U_i \cap U_j, \Ic)
$$
we get
$$
s_j(s_i - s_j) = 0 ,
$$
hence
$$
(1+s_j)(1+s_i-s_j)=1+s_i,
$$
i.e.
$$
1+s_i-s_j = \frac{1+s_i}{1+s_j}.
$$
This finishes the proof. \hfill $\Box$

% ----------------------------------------------------------------------------------------------------------------------------------------------------------------------------------------------------------------------------------------
\subsection{The multiplicity of a rational singularity}

The following theorem is the main result of this section.

\begin{theorem}
\label{theo3} {\rm (\cite{A1} cor. 6)} The multiplicity of the
rational singularity $(X,x)$ is equal to the negative of the
self-intersection-number of the numerical cycle $-(\Znum)^2$.
\end{theorem}

For the proof, we need two lemmas, which are interesting in their own rights.

\begin{lemma}
\label{mznum} {\rm (\cite{A1} thm. 4)} We have

$$ \mm _x \cdot \Oxx = \Oxx(-\Znum) . $$
\end{lemma}

{\bf Proof} (\cite{A1}, \cite{Reid} 4.17):

The inclusion $ \mm _x \cdot \Oxx \subseteq \Oxx(-\Znum)$
is easy \cite{A1}: For $f \in \Gamma(U,\mm _x \cdot \Oxx)$ we can
split the principal divisor $(f)$ in a part $Z$ supported on $E$
and a part $D$, which does not involve any of the $E_i$ at all:
$(f) = Z + D$. Obviously $Z > 0$. We have
$$
(f) \cdot E_i = 0 \text{ and } D \cdot E_i \ge 0 \text{ for } i=1, \dots , n,
$$
since $f$ is regular in a neighbourhood of $E$. Thus $Z \cdot E_i \le 0$ for $i=1, \dots, n$, that is $Z \ge \Znum$. Hence
$f \in \Gamma(U,\Oxx(-\Znum))$.

For the other inclusion we have to show that for each point $p \in
E$ there exists a local section $f$ of $\mm _x \cdot \Oxx$ such that $(f)_{|U} =
{\Znum}_{|U} $ for a neighbourhood $U$ of $p$ (\cite{Reid} 4.17).

Let $X^{\prime}$ be an affine neighbourhood of $x \in X$ and set
$\widetilde{X}^{\prime} := \widetilde{X} \times _X X^{\prime} = \pi^{-1}(X^{\prime})$.

We will write for short ${\II} := \OO_{\widetilde{X}^{\prime}}(-\Znum)$.

We can construct a divisor $A$ on \Znum\ as a sum of local transverse cuts such that $p \not\in \text{Supp} A$ and $\degz A = \degz \II_{|\Znum} $.
The crucial point is, that proposition \ref{picz} implies now $\OO_{\Znum}(A) \cong \II_{|\Znum} $.
Hence there exists a section $s \in H^0(E, \II_{|\Znum})$ which does not vanish at $p$.

To finish the proof, all we have to do is to lift $s$ to a section on $\widetilde{X}^{\prime}$.

From the short exact sequence

$$ 0 \rightarrow \II^{\otimes 2} \rightarrow \II \rightarrow \II_{|\Znum} \rightarrow 0 $$

(obtained from $ 0 \rightarrow \II \rightarrow \OO_{\widetilde{X}^{\prime}} \rightarrow \OO_{\Znum} \rightarrow 0 $ by tensoring with \II)
we get

$$ H^0(\widetilde{X}^{\prime},\II) \rightarrow H^0(E,\II_{|\Znum}) \rightarrow H^1(\widetilde{X}^{\prime},\II^{\otimes 2}). $$

So it suffices to proof $H^1(\widetilde{X}^{\prime},\II^{\otimes 2}) =0$. We will prove more generally:

\subparagraph{Useful fact:} $H^1(\widetilde{X}^{\prime},\II^{\otimes k}) = 0 \text{ for all } k \in {\mathbb N}$.

{\bf Proof of the useful fact}: Since $X^{\prime}$ is affine

$$ H^1(\widetilde{X}^{\prime},\II^{\otimes k}) = H^0(X^{\prime},R^1 \pi_{\ast} \II^{\otimes k}) $$

by \cite{Ha} III.8.5.

The sheaf $R^1 \pi_{\ast} \II^{\otimes k}$ is concentrated in $x$; hence it is enough to prove
$\left( R^1 \pi_{\ast} \II^{\otimes k} \right)_x = 0 $.
By Grothendieck's theorem on formal functions (\ref{app3})

$$ \left( \left( R^1 \pi_{\ast} \II^{\otimes k} \right)_x \right) \widehat{}\, =
\mathop{\varprojlim}_Z H^1(E, {\II^{\otimes k}}_{|Z} )$$

(cf. the proof of theorem \ref{theo1}), so we are left to show $H^1(E,{\II^{\otimes k}}_{|Z})=0$ for all positive cycles $Z$. Again, we can construct a divisor $A$ on $Z$
as a sum of local transverse cuts such that $\text{deg}_Z A = \text{deg}_Z \II^{\otimes k}$.

From the short exact sequence

$$ 0 \rightarrow \OO_Z \rightarrow \OO_Z (A) \rightarrow \OO_A(A) \rightarrow 0 $$

we get

$$ H^1(E, \OO_Z) \rightarrow H^1(E, \OO_Z(A)) \rightarrow H^1(A,\OO_A(A)) . $$

We have $H^1(E, \OO_Z) = 0$ by corollary \ref{H10} and $H^1(A, \OO_A(A) ) = 0$ for dimension reasons \ref{app6},
hence also $ H^1(E, \OO_Z(A)) = 0 $.
\hfill $\blacksquare$

This shows

$$ H^0(\widetilde{X}^{\prime},\OO_{\widetilde{X}^{\prime}}(-\Znum)) \twoheadrightarrow H^0(E,\II_{|\Znum})$$

and we can fetch a preimage
$ s^{\prime} \in H^0(\widetilde{X}^{\prime},\OO_{\widetilde{X}^{\prime}}(-\Znum))$ of $s$. By construction, ${s^{\prime}}_p \neq 0$, thus
${s^{\prime}}_q \neq 0$ for all $q \in U$ for some neighbourhood $U$ of $p$. Or put differently,

$$ (s^{\prime})_{|U} = {\Znum}_{|U} . $$

We have already seen $\pi_{\ast} \Oxx = \Ox$; therefore $s^{\prime}$ gives rise to a section in $\Gamma(\pi(U), \Ox)$,
and thus in
$\Gamma(\pi(U), \mm _x)$, since $\Znum \ge E$. \hfill $\Box$

\begin{lemma}
\label{gen} {\rm (\cite{Reid} 4.18)} The ring $\bigoplus_{k \ge 0}
H^0(E,\II^{\otimes k})$ is generated in degree $1$, where
$\II:=\OO_{\Znum}(-\Znum)$.
\end{lemma}

{\bf Proof}:
We can use local transverse cuts to construct a divisor $A$ on \Znum\ with $\OO_{\Znum}(A) \cong \II$. Therefore there exists a global section
$s_0 \in H^0(E, \II)$ whose divisor of zeros is precisely $A$. Since we had a lot of freedom in choosing $A$, we see that the linear system
$|A| = |\II |$ is basepoint-free.
Thus we can choose a $s \in H^0(E,\II)$ such that $s$ provides a local base at every point $q \in A$.
We can use $s^{k-1}$ to identify
$\II^{\otimes k-1} \otimes \OO_A \cong \OO _A$.
The short exact sequence

$$ 0 \rightarrow \OO_{\Znum} \rightarrow \II \rightarrow \OO _A \rightarrow 0 $$

yields (corollary \ref{H10})

\begin{equation}
\label{ses1}
0 \rightarrow H^0(E,\OO _{\Znum}) \rightarrow H^0(E,\II) \rightarrow H^0(A,\OO _A) \rightarrow 0
\end{equation}

and (\ref{app6})

$$ H^1(E, \II)=0 . $$

Let $s_1, \dots, s_d \in H^0(E, \II)$ map to a basis of $H^0(A,\OO _A)$, $d:=h^0(A,\OO _A)$.
Our lemma will follow from the following claim:

\begin{eqnarray*}
& H^1(E,\II ^{\otimes k})=0, \\
& H^0(E,\II ^{\otimes k})=\text{span}_{\CC}\{{s_0}^k, {s_0}^{k-l}s^{l-1}s_i : 1 \le l \le k, 1 \le i \le d \}
\end{eqnarray*}

Note that the sections
$ {s_0}^{k-l}s^{l-1}s_i ,  l = 1, \dots,  k, i = 1, \dots, d$
and ${s_0}^k$
are linearly independent over \CC.

For $k=1$ our claim follows from (\ref{ses1}) and $h^0(E,\OO_{\Znum})=1$:
$$
h^0(E,\II)=h^0(E,\OO _{\Znum}) + h^0(A,\OO _A) = 1 + d.
$$
For the induction step we argue similarly using
$$
\begin{array}{ccccccc}
0 & \rightarrow & H^0(E,\II ^{\otimes k-1}) & \rightarrow & H^0(E,\II ^{\otimes k}) & \rightarrow & H^0(A,\OO_A) \\
& \rightarrow & 0 & \rightarrow & H^1(E, \II ^{\otimes k}) & \rightarrow & 0
\end{array}
$$
to get
$$
h^0(E,\II ^{\otimes k}) = h^0(E,\II ^{\otimes k-1}) + h^0(A,\OO _A) = (k-1)d+1+d
$$
and
$$
H^1(E,\II ^{\otimes k})=0 .
$$
\hfill $\Box$

{\bf Proof of the theorem} (\cite{Reid} 4.18): Let $X^{\prime}=\text{Spec }R$ be
an affine neighbourhood of $x \in X$ and set
$\widetilde{X}^{\prime}:=\widetilde{X} \times _X X^{\prime}$. By
lemma \ref{mznum}, we know
\begin{equation}
\label{k1}
H^0(\XXX, \mm _x \cdot \Oxx) = H^0(\XXX, \Oxx(-\Znum)) .
\end{equation}

We want to generalize (\ref{k1}) to

\begin{equation}
\label{kk}
H^0(\XXX, {\mm _x }^k \cdot \Oxx) = H^0(\XXX, \Oxx(-k\Znum)) .
\end{equation}

With (\ref{kk}) the proof of our assertion is straightforward:
Since
$$
H^1(\XXX, \Oxx(-(k+1)\Znum))=0
$$
by corollary \ref{H10}, we have
$$
\frac{H^0(\XXX,{\mm _x} ^k \cdot \Oxx)}{H^0(\XXX,{\mm _x} ^{k+1} \cdot \Oxx)} =
\frac{H^0(\XXX,\Oxx(-k\Znum))}{H^0(\XXX,\Oxx(-(k+1)\Znum))} =
H^0(E, \OO _{\Znum} (-k\Znum)) .
$$
The Riemann-Roch theorem for curves tells us
\begin{eqnarray*}
h^0(E,\OO _{\Znum}(-k\Znum)) & = & 1 - p_a(\Znum) + \text{deg } \OO _{\Znum} (-k\Znum) \\
& = & 1 - k (\Znum)^2,
\end{eqnarray*}
that is the leading coefficient of the Hilbert-Samuel polynomial of $(\OO _{X^{\prime} , x}, m_x)$ is $-(\Znum)^2$.

We will prove (\ref{kk}) by induction, so assume (\ref{kk}) holds for $k < l$.
Clearly (\ref{k1}) implies the inclusion

$$ H^0(\XXX, {\mm _x} ^l \cdot \Oxx) \subseteq H^0(\XXX, \Oxx(-l\Znum)). $$

We want to show surjectivity. We take a $g \in H^0(\XXX, \Oxx(-l\Znum))$ and restrict it to
$ \bar{g} \in H^0(E, \OO _{\Znum} (-l\Znum)). $ By lemma \ref{gen} we have a surjection

$$ H^0(E,\OO _{\Znum} (-\Znum)) \otimes H^0(E,\OO _{\Znum} (-(l-1)\Znum)) \twoheadrightarrow H^0(E,\OO _{\Znum} (-l\Znum)), $$

i.e. we can write $\bar{g}$ in the form $\bar{g} = \sum\limits_{j=1}^m \bar{x}_j \bar{y}_j $ with
\begin{eqnarray*}
&& \bar{x}_j \in H^0(E,\OO _{\Znum} (-\Znum))  \text{ and }\\
&& \bar{y}_j \in H^0(E,\OO _{\Znum} (-(l-1)\Znum)) \text{ for } j = 1,\dots, m .
\end{eqnarray*}
Lifting $\bar{x}_j$ and $\bar{y}_j$ to sections on \XXX
\begin{eqnarray*}
&& x_j \in H^0(\XXX,\Oxx (-\Znum)) \text{ and } \\
&& y_j \in H^0(\XXX,\Oxx (-(l-1)\Znum)) \text{ for } j = 1,\dots, m
\end{eqnarray*}
gives for
$f_2:=\sum\limits_{j=1}^m x_j y_j$ by induction hypothesis
\begin{eqnarray*}
& f_2 \in H^0(\XXX,{\mm _x} ^l \cdot \Oxx) \text{ and }\\
& g-f_2 \in H^0(\XXX, \Oxx (-(l+1)\Znum)).
\end{eqnarray*}
Continuing in this fashion gives
\begin{eqnarray*}
& f_2, \dots , f_p \in H^0(\XXX,{\mm _x} ^l \cdot \Oxx) \text{ and }\\
& g-f_2- \dots -f_p \in H^0(\XXX, \Oxx (-(l+p-1)\Znum))  .
\end{eqnarray*}

Hence it suffices to prove
$$
H^0(\XXX, \Oxx (-p\Znum)) \subseteq H^0(\XXX,{\mm _x} ^l \cdot \Oxx) \text{ for } p \gg 0 .
$$
The point is that
$$
\bigoplus_{k \ge 0} H^0(\XXX, \Oxx(-k\Znum))
$$
is a finitely generated $R$-algebra. Assuming this, let $M$ be the maximal degree in a fixed set of generators. For $p > l M$, each element of
$H^0(\XXX, \Oxx(-p\Znum))$ is a sum of products of at least $l$ generators, thus
$$
H^0(\XXX, \Oxx (-p\Znum)) \subseteq H^0(\XXX,{\mm _x} ^l \cdot \Oxx) \text{ for } p > l M  .
$$
For the proof of the assumption, note that the complete linear system $| \OO _{\XX ^{\prime}}(-\Znum) |$ is free:
By the useful fact, $H^1(\XX ^{\prime} , \OO _{\XX ^{\prime}}(-2\Znum))=0$, i.e. we have a surjection
$$
H^0(\XX ^{\prime} , \OO _{\XX ^{\prime}}(-\Znum)) \twoheadrightarrow H^0(E,\OO _{\Znum}(-\Znum)) .
$$
We have seen in the proof of lemma \ref{mznum} that  $| \OO _{\Znum}(-\Znum) |$ is free and hence so is
$| \OO _{\XX ^{\prime}}(-\Znum) |$.
Thus we have a well-defined morphism
$$
\phi_{|- \Znum |} : \XX ^{\prime} \rightarrow X^{\prime} \times {\mathbb P}^N = {\mathbb P}_{X^{\prime}}^N .
$$
We denote its image by $Y:=\text{im }  \phi_{|- \Znum |} $. We want to show that $Y$ is closed. Since $\XX \rightarrow X$ is proper,
so is $\XX^{\prime} \rightarrow X^{\prime}$ (\cite{Ha} II.4.8.(c)). Because of the separatedness of
${\mathbb P}_{X^{\prime}}^N \rightarrow X^{\prime}$ (\cite{Ha} II.4.9), we see that
$\phi_{|- \Znum |} : \XX ^{\prime} \rightarrow  {\mathbb P}_{X^{\prime}}^N $ is proper (\cite{Ha} II.4.8.(e)), in particular
$\phi_{|- \Znum |}$ is closed. Hence $Y$ is closed.

The pullback under $ \phi_{|- \Znum |}  $
of the relatively ample line bundle $\OO(1) := \OO _{X^{\prime}} \otimes _{\CC} \OO_{{\mathbb P}^N}(1)$ is by
definition of $\phi_{|- \Znum |}$ simply $\OO_{\XX^{\prime}} (-\Znum)$ (\cite{Reid} 4.18).

Thus it suffices to show that
$$
S^{\prime}(Y) := \bigoplus_{k \ge 0} H^0(Y,\OO _Y  (k))
$$
is finitely generated as a $R$-algebra. The homogenous coordinate ring
$S(Y)=A[x_0, \dots , x_N]/I(Y)$
of $Y$ is certainly a finitely generated $R$-algebra.
By \cite{Ha} Ex. II.5.9, there exists a natural graded morphism
$$
S(Y) \rightarrow S^{\prime}(Y) ,
$$
which is an isomorphism in high degrees, i.e.
$$
S(Y)^d \stackrel{\sim}{\longrightarrow}  S^{\prime}(Y)^d \text{ for } d \gg 0,
$$
and we are done. \hfill $\Box$

\begin{corollary}
\label{znum2}
For a rational double point $(X,x)$, we have
$(\Znum)^2 = -2$.
\end{corollary}

% -------------------------------------------------------------------------------------------------------------------------------------------------------------------------------------------------------------------------------------------
\section{The geometry of the exceptional set $E$ of a re\-solution of a rational double point}

Once the hard work has been done in proving theorems \ref{theo1} and \ref{theo3}, it is now easy to say explicitly
what configurations can arise for $E$, if $(X,x)$ is a rational double point.

From now on, we will assume that \res\ is a good resolution of a rational double point $(X,x)$ and $E$ its exceptional set.

By proposition \ref{negdef}, we have $E_i^2 \le -1$ for $i=1, \dots ,
n$. If $E_{i_0}^2 = -1 $ for some $i_0$, then $E_{i_0} \cong \PP$
can be contracted by Castelnuovo's criterion \ref{app4} to give a
resolution $\pi^{\prime}: \XXX \rightarrow X$ with fewer $E_i$.
(In general, the resolution $\pi^{\prime} : \XXX \rightarrow X$
needs not to be good anymore, since the condition \ref{goodcond2} in the definition of a good resolution
might be violated. However, it is a simple consequence of the
following theorem \ref{theo4} and, again, lemma \ref{negdef}, that
this cannot happen in our case. Note that we do not use this condition
\ref{goodcond2} in the proof of \ref{theo4}.)

Therefore, we can assume $E_i^2 \le -2$ for $i=1, \dots, n$ without loss of generality.
\begin{theorem}
\label{theo4} {\rm (\cite{Df})} The $E_i$ have
self-intersection-number $-2$.
\end{theorem}

{\bf Proof}: Let $K$ be a canonical divisor on $\widetilde{X}$
(\cite{Ha} V.1.4.4). The adjunction formula \ref{app2} tells us

\begin{equation}
\label{adjdf}
-E_i \cdot K = E_i^2 + 2
\end{equation}

and thus

$$ E_i \cdot K \ge 0 . $$

We apply the adjunction formula \ref{app2} for \Znum

$$ 2 p_a(\Znum)-2 = (\Znum)^2 + \Znum \cdot K, $$

and get by the corollaries \ref{znum0} and \ref{znum2}

$$ 0 = \Znum \cdot K = \sum\limits_{i=1}^n r_i(E_i \cdot K) \ge 0,
\text{ i.e. } E_i \cdot K =0. $$

Using (\ref{adjdf}) again, we see $E_i^2=-2$. \hfill $\Box$

\vspace{3mm} We define the {\em Dynkin diagram} of the resolution
\res\ to be the weighted dual graph $\Gamma$ associated to $E$:
The vertices $e_i$ of $\Gamma$ correspond to the $E_i$. Whenever $E_i$
and $E_j$ intersect for $i \neq j$, the corresponding vertices are
joined by an edge. Finally, we associate to every vertex $e_i$ of
$\Gamma$ the self-intersection-number $E_i^2$.

Every weighted graph $\Gamma$ defines a bilinear form $\langle
\cdot,\cdot \rangle$ on the free module with the vertices $e_i,
i=1, \dots, n$ of $\Gamma$ as basis in the following way: We take
\begin{eqnarray*}
& \langle e_i,e_i \rangle := \text{the weight of }e_i \text{ and }
\\ & \langle e_i,e_j \rangle := \text{number of edges joining }e_i
\text{ and }e_j.
\end{eqnarray*}

The bilinear form of the Dynkin diagram of a resolution is obviously given by the matrix
$(E_i \cdot E_j)_{i,j=1,\dots,n}$ and hence negative definite by proposition
\ref{negdef}. This puts very strong restrictions on the possible Dynkin diagrams $\Gamma$.

\begin{proposition} {\rm (\cite{Df})}
Let $\Gamma$ be a connected graph weighted by $-2$ whose
associated bilinear form is negative definite. Then $\Gamma$ is a
$T$-tree $T_{p,q,r}$

\begin{center}
\Ttree
\end{center}

with ${1 \over p} + {1 \over q} + {1 \over r}
> 1$.
\end{proposition}

{\bf Proof}: Every connected subgraph $\Gamma^{\prime}$ of
$\Gamma$ satisfies the hypothesis as well, hence can be
\begin{itemize}
\item neither a loop

\begin{center}
\cyclus,
\end{center}

since then $(e_1 + \dots + e_n)^2=0$ contradicting the negative
definiteness condition
\item nor of the form

\begin{center}
\cross,
\end{center}

since then $(2e_1 + \dots + 2e_n + f_1 + \dots + f_4)^2=0$.
\end{itemize}

Thus $\Gamma$ must be of the form $T_{p,q,r}$.

The condition ${1 \over p} + {1 \over q} + {1 \over r} > 1$ follows by an elementary argument: With respect to the standard basis given by the vertices
of $\Gamma$, the associated bilinear form is expressed by the matrix

$$
\left(
\begin{array}{ccccccccccc}
-2 & 1 &&&&&&&&& \\
1 & -2 & 1 &&&&&&&& \\
& 1 & -2 &&&& \mathop{\fbox{1}}\limits_{p,p+q} &&&& \\
&&& \ddots &&&&&&& \\
&&&& -2 & 1 &&&&& \\
&&&& 1 & -2 & \mathop{\fbox{0}}\limits_{p+q-1,p+q} &&&& \\
&&& \mathop{\fbox{1}}\limits_{p+q,p} && \mathop{\fbox{0}}\limits_{p+q,p+q-1} & -2 & 1 && \\
&&&&&& 1 & -2 & 1 && \\
&&&&&&& 1 & -2 && \\
&&&&&&&&& \ddots & 1 \\
&&&&&&&&& 1 & -2
\end{array}
\right).
$$

But, up to congruence, this is equal to

$$
\left(
\begin{array}{cccccccccc}
-\frac{2}{1} &&&&&&&&& \\
& -\frac{3}{2} &&&&&&&& \\
&& \ddots &&&&&&& \\
&&& -\frac{p+1}{p} & \mathop{\fbox{1}}\limits_{p,p+1} &&& \mathop{\fbox{1}}\limits_{p,p+q} && \\
&&& \mathop{\fbox{1}}\limits_{p+1,p} & -\frac{q}{q-1} &&&&& \\
&&&&& \ddots &&&& \\
&&&&&& -\frac{2}{1} &&& \\
&&& \mathop{\fbox{1}}\limits_{p+q,p} &&&& -\frac{r}{r-1} && \\
&&&&&&&& \ddots & \\
&&&&&&&&& -\frac{2}{1}
\end{array}
\right).
$$

Now, this matrix is congruent  to a diagonal matrix with negative
main diagonal entries, except maybe a single one
$1-p^{-1}-q^{-1}-r^{-1}$. \hfill $\Box$

\begin{corollary}
The Dynkin diagram associated to a rational double point $(X,x)$
must be one of the following diagrams
{\rm
\begin{center}
$A_n$ \raisebox{-2mm}{\An} ($n$ vertices), \\
$D_n$ \raisebox{-11mm}{\Dn}($n$ vertices), \\
$E_6$ \raisebox{-11mm}{\Esechs},\hspace{5mm} $E_7$ \raisebox{-11mm}{\Esieben} and \\
$E_8$ \raisebox{-11mm}{\Eacht}.
\end{center}
}
\end{corollary}

We will see in section \ref{C2G} that all Dynkin diagrams actually occur. We say that a rational double point is of type $A_n$, $D_n$ or $E_n$
according to its Dynkin diagram.

% -------------------------------------------------------------------------------------------------------------------------------------------------------------------------------------------------------------------------------------------------
\section{Example: The singularities $\CC/G$ for finite $G \subset \text{SL}(2,\CC)$}
\label{example1}
After a bit of the theory of rational double points has been presented, we want to study an example.

% ----------------------------------------------------------------------------------------------------------------------------------------------------------------------------------------------------------------------------------------
\subsection{Conjugacy classes of finite subgroups of $\text{SL}(2,\CC)$}

As a preliminary, we recall briefly the classification of
conjugacy classes of finite subgroups of $\text{SL}(2,\CC)$. We
consider first $\text{SO}(3,{\mathbb R})$. Up to conjugacy, the
finite subgroups of $\text{SO}(3,{\mathbb R})$ are the rotational
symmetry groups of
\begin{itemize}
\item a pyramid (giving the {\em cyclic subgroups} $C_n$)

\begin{center}
\includegraphics[scale=0.1]{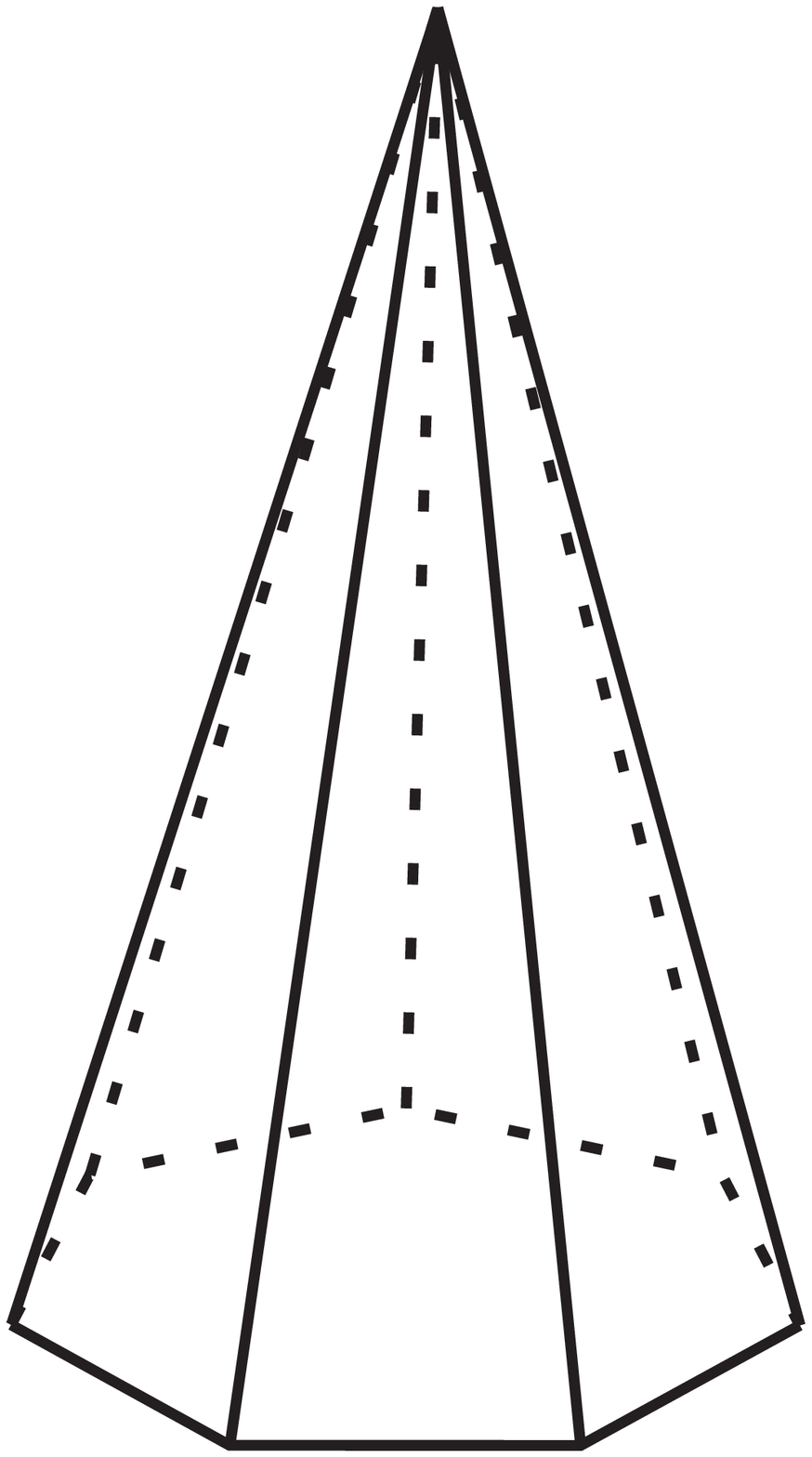},
\end{center}

\item an orange (corresponding to the {\em dihedral subgroups} $D_n$)

\begin{center}
\includegraphics[scale=0.3]{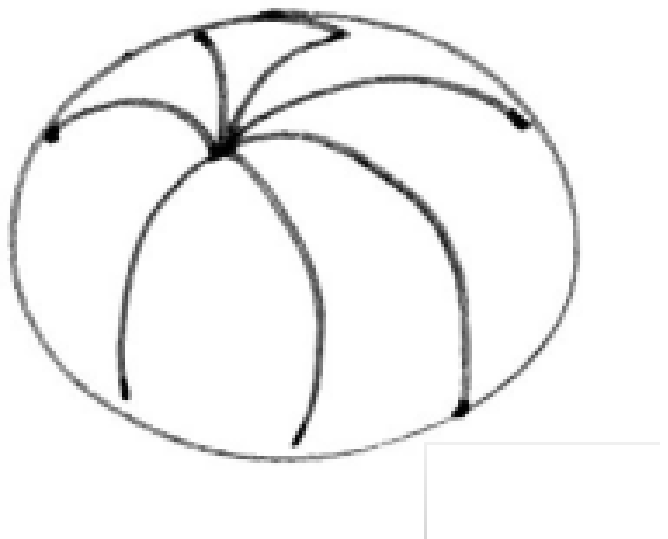},
\end{center}

\item and the Platonic solids, which give
\begin{itemize}
\item the {\em tetrahedral subgroup} $T=A_4$

\begin{center}
\includegraphics[scale=0.35]{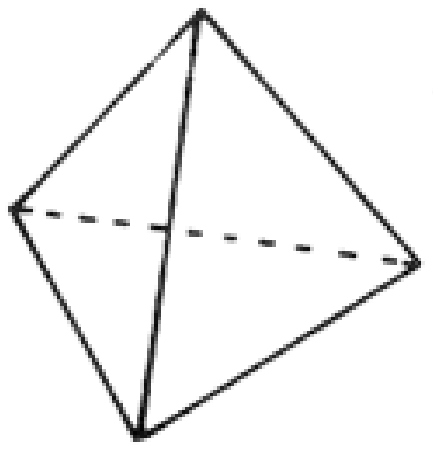},
\end{center}

\item the {\em octahedral subgroup} $O=S_4$

\begin{center}
\includegraphics[scale=0.35]{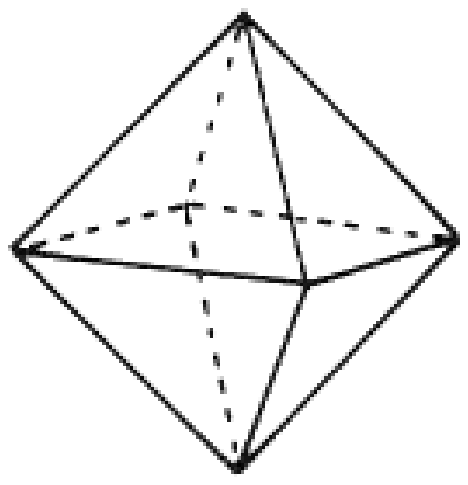},
\end{center}

\item and the {\em icosahedral subgroup} $I=A_5$

\begin{center}
\includegraphics[scale=0.35]{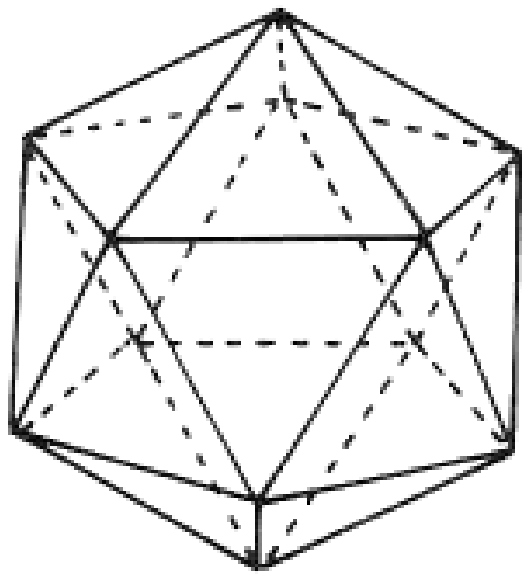},
\end{center}

\end{itemize}
\end{itemize}
respectively \cite{Sl}.

If we identify $S^2 \cong \PP$, we get an inclusion of the group of isometries of $\PP$ (with respect to the usual metric) into the group of
conformal transformations

$$ \text{SO}(3,{\mathbb R}) \subset {\mathbb P}\text{GL}(2,\CC) . $$

Under the double cover

$$ \rho: \text{SL}(2,\CC) \twoheadrightarrow {\mathbb P}\text{GL}(2,\CC) = \text{SL}(2,\CC) / \{\pm 1\}  $$

this inclusion corresponds to

$$ \rho^{-1}(\text{SO}(3,{\mathbb R})) = \text{SU}(2,\CC). $$

Since for any finite subgroup $G$ of \text{SL}(2,\CC) we can find
a $G$-invariant Hermitian metric by averaging an arbitrary one,
every finite subgroup $G$ of $\text{SL}(2,\CC)$ is conjugated to a
subgroup of $\text{SU}(2,\CC)$. Hence it corresponds to a finite
subgroup of $\text{SO}(3,{\mathbb R})$, unless it is a cyclic
group of odd order. Thus we have derived the following
classification of the conjugacy classes of finite subgroups of
$\text{SL}(2,\CC)$ (\cite{Lm} II \S 1):
\begin{itemize}
\item the cyclic subgroup of order $n$ $C_n$,
\item the binary dihedral subgroups $\widetilde{D}_n = \rho^{-1}(D_n)$,
\item the binary tetrahedral, octahedral and icosahedral subgroup $\widetilde{T}=\rho^{-1}(T)$, $\widetilde{O}=\rho^{-1}(O)$
and $\widetilde{I}=\rho^{-1}(I)$ respectively.
\end{itemize}

% ----------------------------------------------------------------------------------------------------------------------------------------------------------------------------------------------------------------------------------------
\subsection{The singularities $\CC^2 / G$}
\label{C2G}

Now let $G$ be any of these subgroups; the affine orbit variety

$$ \CC^2 / G = \text{Spec }\CC[ x_1, x_2]^G $$

has an isolated singularity at the origin. The singularities
obtained in this fashion are all rational double points \cite{Df}.
It is a result from classical invariant theory that these
singularities embed in codimension one

$$ \CC^2 / G = \text{Spec }\CC[x,y,z] / (f) \text{ where } f \in \CC[x,y,z] . $$

See for example Klein's influential book \cite{Klein}, or also \cite{DuVal3} 5.39. For a modern treatment, we refer to \cite{DPT}, p. 5.
The following table \ref{rdptable} contains the basic information about these singularities.

\begin{table}[h]
\caption{The singularities $\CC^2 / G = \text{Spec }\CC[x,y,z] /
(f) $ for finite $G \subset \text{SL}(2,\CC)$} \label{rdptable}
\begin{minipage}{\linewidth}
\begin{tabular}{|c|c|c|c|c|}
\hline
$G$ & $f$ & Name & Dynkin diagram &
\Znum\footnote{The number above a vertex denotes the multiplicity of the corresponding
projective line $E_i \cong \PP$ in the numerical cycle \Znum.} \\
\hline \hline
&&& \An & \Anznum \\
$C_n$  & $x^n+y^2+z^2$ & $A_{n-1}$ &
($n-1$ vertices) & \\
\hline
&&& \Dn & \Dnznum \\
$\widetilde{D}_n$& $x^{n+1}+xy^2+z^2$ & $D_{n+2}$ &
($n+2$ vertices) & \\
\hline
$\widetilde{T}$ & $x^2+y^3+z^4$ & $E_6$ &
\Esechs & \Esechsznum \\
\hline
$\widetilde{O}$ & $x^2+y^3+yz^3$ & $E_7$ &
\Esieben & \Esiebenznum \\
\hline
$\widetilde{I}$ & $x^2+y^3+z^5$ & $E_8$ &
\Eacht & \Eachtznum \\
\hline
\end{tabular}
\end{minipage}
\end{table}

% ----------------------------------------------------------------------------------------------------------------------------------------------------------------------------------------------------------------------------------------
\subsection{The icosahedral case $G=\widetilde{I}$}
\label{icosahedral}

We will sketch the proof of the assertions made so far in the special case $G=\widetilde{I}$. We see that

$$ X:= \CC^2 / \widetilde{I} = \text{Spec }\CC [ x, y, z ] / (x^2 + y^3 + z^5) $$

has a singularity at $x_0:=(0,0,0)$, which must be a double point, since
$$\text{rank}_{\CC}\frac{(x,y,z)^k+(x^2+y^3+z^5)}{(x,y,z)^{k+1}+(x^2+y^3+z^5)} \le 2k+1 .$$

We want to show that $(X,x_0)$ is a {\bf rational} double point. Indeed, we will show how to resolve the singularity $(X,x_0)$ (\cite{Lm} IV \S 9,
\cite{DPT} p. 15). Blowing up $\CC^2$ instead of $\CC^2 / \widetilde{I}$ will reveal the relevant information more clearly. The blow-up of $\CC^2$
at the origin is known to be ${\mathbf v}(\OO _{\PP} (-1))$ (\cite{Ha} V.3.1);
here ${\mathbf v}({\mathcal F})$ denotes the vector bundle determined by the locally trivial
sheaf $\mathcal F$. The quotient of ${\mathbf v}(\OO _{\PP} (-1))$ by $\{\pm 1\}$ is ${\mathbf v}(\OO _{\PP} (-2))$:

\begin{eqnarray*}
& {\mathbf v}(\OO _{\PP} (-1)) / \{\pm 1\} = \{ (w_1:w_2;z_1,z_2)
\in \PP \times \CC ^2 : w_1 z_2 = w_2 z_1 \} / \{\pm 1\} \cong \\
& \{ (w_1:w_2;z_1,z_2) \in \PP \times \CC ^2 : w_1^2 z_2 = w_2^2
z_1 \} =  {\mathbf v}(\OO _{\PP} (-2)) \\ & (w_1:w_2;z_1,z_2)
\mapsto (w_1:w_2;z_1^2,z_2^2) .
\end{eqnarray*}
Note that the $\widetilde{I}$-action on $\CC^2$ lifts to an action on ${\mathbf v}(\OO _{\PP} (-1))$.
Thus we have the following commutative diagram

$$
\begin{CD}
{\mathbf v}(\OO _{\PP}(-2)) \cong {\mathbf v}(\OO _{\PP}(-1)) / \{ \pm 1 \}  @<<< {\mathbf v}(\OO _{\PP}(-1)) @>\sigma>> \CC^2 \\
@VVV @VVV @VVV \\
{\mathbf v}(\OO _{\PP}(-2)) / I @<\sim<< {\mathbf v}(\OO _{\PP}(-1)) / \widetilde{I} @>\bar{\sigma}>> \CC^2 / \widetilde{I}
\end{CD}.
$$

In particular, $\bar{\sigma}^{-1}(x_0)$ is a copy of \PP\ with self-intersection-number $-2$. A precise analysis of the $I$-action on
${\mathbf v}(\OO _{\PP} (-2))$ shows, that $I$ acts on the zero-section $S^2 \cong \PP \subseteq {\mathbf v}(\OO _{\PP} (-2))$ in the usual way
as rotations which leave  an inscribed icosahedron invariant. Furthermore the action of $I$ on $T^{\ast}\PP \cong {\mathbf v}(\OO _{\PP} (-2))$
is simply the cotangent action induced by the action of $I$ on $\PP$ (\cite{Lm} IV \S 7). Now, the group $I$  is acting free on $S^2$, except on three
exceptional orbits, which consist of the vertices, the mid-edge-points and the mid-face-points of the inscribed icosahedron, respectively.
Moreover, we see that these three orbits are also the only exceptional orbits for the action of $I$ on ${\mathbf v}(\OO _{\PP} (-2))$.
Therefore, the quotient variety
${\mathbf v}(\OO _{\PP} (-2)) / I $
is smooth except at the three points corresponding to these orbits.
An explicit calculation using local coordinates shows that these three singular points are cyclic quotient singularities of type $(5,4)$, $(3,2)$ and
$(2,1)$, respectively (\cite{DPT} p. 17), this notion being defined as follows: A {\em cyclic quotient singularity of type} $(n,q)$ is a singularity, which
is isomorphic to $\CC^2 / {{\mathbb \mu}_{n,q}}$ where ${{\mathbb \mu}_{n,q}}$ is the cyclic group generated by
$\left(
\begin{array}{cc}
\xi & 0 \\
0 & \xi^q
\end{array}
\right)$
for a $n$-th root of unity $\xi$.
Note that the numbers $5$, $3$, $2$ are the ramification indices at the ramification points of the map

$$ S^2 \twoheadrightarrow S^2 / I $$

corresponding to the exceptional orbits.

A cyclic quotient singularity of type $(n,q)$ can be resolved by the Hirzebruch-Jung algorithm using successive blow-ups of points (\cite{Fulton} 2.6).
The exceptional set of a resolution of a cyclic quotient singularity obtained in this way is a bunch of rational curves; the associated Dynkin
diagram is of the form

\begin{center}
\HiJu ,
\end{center}

where the $b_i$'s are calculated by a modified Euclidean algorithm

$$ {n \over q} = b_1 - {1 \over {b_2 - {1 \over \ddots - {1 \over b_k }}}} .$$

Applying the Hirzebruch-Jung algorithm three times for the three singular points we got, we obtain a resolution
\res\ with associated Dynkin diagram

\begin{center}
\Eachtlab
\end{center}

The numerical divisor \Znum\ is easily verified to be

$$ \Znum = 2E_1 + 4E_2 + 6E_3 + 3E_4 + 5E_5 + 4E_6 + + 3E_7 + 2E_8 . $$

A straightforward computation using equation (\ref{lichtenbaum}) shows

$$ p_a(\Znum) = 0, $$

i.e. the singularity $(X,x_0)$ is rational by the criterion of theorem \ref{theo2}.

% --------------------------------------------------------------------------------------------------------------------------------------------------------------------------------------------------------------------------------------------------
\section{Changing to the complex analytic category}
\label{complexana}

The examples introduced in the preceding section already exhaust
all possibilities of rational double points up to isomorphism in
the complex analytic category. Of course, such a statement cannot
be true in the category of algebraic varieties, since rational
double points can live on the various kinds of surfaces and the
birationality class of the surface is encoded locally due to the
coarseness of the Zariski topology.

We give the following (simplified) definition of the {\em complex analytic category}:
\begin{itemize}
\item Its objects are called {\em complex analytic spaces} and can be constructed as follows. Let $U$ be a simply connected open subset of $\CC^n$,
$\OO _U$ the sheaf of complex analytic functions on $U$, and $\II _X$ a sheaf of ideals on $(U,\OO _U)$. Denote by $X \subseteq U$ the zeroset
of $\II _X$ equipped with the {\bf standard topology} and set $\OO _X := \OO _U / \II _X$. The pair $(X,\OO _X)$ is then a complex analytic space.
(For the general definition one allows such simple building blocks to be glued together as in the
definition of schemes.)
\item a {\em morphism of complex analytic spaces} $f: (X,\OO _X) \rightarrow (X^{\prime},\OO _{X^{\prime}})$ is a continuous map
$f: X \rightarrow X^{\prime}$ such that $f^{\ast}: \OO_{X^{\prime}} \rightarrow \OO _X$ is well-defined.
\end{itemize}
To a great extent, the complex analytic category is similar to the category of algebraic varieties: For example, stalks $\OO _{X,x}$ are
noetherian local rings and for reduced complex analytic spaces $(X,\OO _X)$ R\"{u}ckert's Nullstellensatz holds (\cite{Lm} III \S 8).

% ----------------------------------------------------------------------------------------------------------------------------------------------------------------------------------------------------------------------------------------
\section{Tautness of rational double points}

% ----------------------------------------------------------------------------------------------------------------------------------------------------------------------------------------------------------------------------------------
\subsection{Definition and theorem}

To pick up the question of classifying rational double points in the complex analytic category, we introduce the notion of tautness.

Let $(X,x)$ be a two-dimensional normal singularity with a good resolution,
whose exceptional set is a bunch of rational curves \PP\, and let
$\Gamma$ be its Dynkin diagram.

We say $(X,x)$ is {\em taut}, if up to analytic isomorphism,
$(X,x)$ is the unique such singularity, that has a good resolution with a bunch of rational curves \PP\ as
exceptional set and $\Gamma$ as its Dynkin diagram \cite{B3}.

We have the following theorem.

\begin{theorem}
\label{rdptaut} The singularities listed in table \ref{rdptable}
are taut.
\end{theorem}

This gives us a complete classification of rational double points up to analytic isomorphism. There are several proofs available for theorem \ref{rdptaut}.

% ----------------------------------------------------------------------------------------------------------------------------------------------------------------------------------------------------------------------------------------
\subsection{Tjurina's proof}

Maybe the most natural one is the proof by Tjurina \cite{Tj}:
Suppose there were another singularity $(X^{\prime},x^{\prime})$ with an exceptional set consisting only of rational curves and
the same Dynkin diagram, i.e. with an isomorphic exceptional variety $(E^{\prime},\OO _{E^{\prime}}) \cong (E,\OO _E)$. Then a sufficient condition for
the existence of an isomorphism of neighbourhoods of $E$ and $E^{\prime}$ is by \cite{Gr} Thm. 3, that $(E,\OO _{n E})$ and
$(E^{\prime},\OO _{n E^{\prime}})$ are isomorphic for $n$ large enough. The proof given in \cite{Tj} proceeds by induction. Assuming
we are given some isomorphism of $(E,\OO _Z)$ and $(E^{\prime},\OO _{Z^{\prime}})$ (here $Z$ and $Z^{\prime}$ are exceptional divisors supported on
$E=\bigcup_{i=1}^n E_i$, $E^{\prime}=\bigcup_{i=1}^n E_i^{\prime}$, respectively), then this can be extended to an isomorphism of $(E,\OO_{Z+E_i})$ and
$(E^{\prime},\OO _{Z^{\prime} + E_i^{\prime}})$, unless some
obstruction occurs, which lies in some cohomology group \cite{Gr}.
Grauert argues, that if all these cohomology groups vanish, then the singularity in question is taut. In general, these cohomology groups
do not vanish and Tjurina's proof is more subtle. Essentially, he shows that the cohomology groups are too small to put obstructions on the lifting
of every possible isomorphism of $(E,\OO _Z)$ and $(E^{\prime},\OO _{Z^{\prime}})$.

% ----------------------------------------------------------------------------------------------------------------------------------------------------------------------------------------------------------------------------------------
\subsection{Brieskorn's first proof}

For the sake of historical correctness, we mention that the first proof of theorem \ref{rdptaut}
was given by Brieskorn (\cite{B4} Satz 1). He showed that
rational double points can be resolved by blowing up points alone, that is,
it is not necessary to normalize or blow up curves. Such singularities
are called {\em absolutely isolated} and were studied by Kirby (\cite{Kirby} 2.6, 2.7),
who gave a classification of absolutely isolated double points:
they are precisely those listed in table \ref{rdptable}.

% ----------------------------------------------------------------------------------------------------------------------------------------------------------------------------------------------------------------------------------------
\subsection{Brieskorn's second proof}
\label{bries}

However, we want to sketch another proof of theorem \ref{rdptaut}, which, also due to Brieskorn \cite{B3}, is of compelling beauty and combines
ideas from different fields of mathematics:

The {\em local fundamental group} of $(X,x)$ is defined as

$$ \pi _{X,x} := \mathop{\varprojlim}_U \pi_1(U \setminus \{x\} )$$

where $U$ runs over all neighbourhoods of $x \in X$ (\cite{B3} \S 2). Equivalently, we can calculate $\pi _{X,x}$ as

$$ \pi _{X,x} \cong \mathop{\varprojlim}_{\widetilde{U}} \pi_1(\widetilde{U} \setminus \{x\} )$$

where the limit is now taken over all neighbourhoods $\widetilde{U}$ of $E \subset \widetilde{X}$. To actually compute $\pi_{X,x}$ it is sufficient
to work out $\pi_1$ for a {\em good} neighbourhood $U$. According to \cite{Prill}, a neighbourhood $U$ of $x \in X$ is called {\em good}, if there exists
a neighbourhood basis $\{ U_i \}$ for $x$ such that $U_i \setminus \{ x \}$ is a deformation retract of $U \setminus \{ x \}$ for all $i$.
Such a good neighbourhood has the homotopy type of a tubular neighbourhood $M$ of $E$ in the sense of Mumford \cite{Mf}.
Intuitively spoken, a tubular neighbourhood is a levelset of the potential distribution due to a uniform charge on $E$. Mumford studied these tubular
neighbourhoods $M$ and showed that they are built out of standard pieces $S^1 \times S^1 \times [ 0,1]$ "plumbed" together in a certain fashion
determined by $(E_i \cdot E_j)_{i,j=1,\dots,n}$.
This description and the Seifert-van-Kampen theorem enables him to give a presentation for $\pi_1 (M)$ in terms of generators and relations. (The ideas
of his proof can also be found in \cite{Lm} IV \S\S 10 - 14 on plumbed surfaces.)
It turns out that for the intersection matrices $(E_i \cdot E_j)_{i,j=1,\dots,n}$
of the resolutions of rational double points this group $\pi _{X,x} = \pi_1 (M)$ is finite.

\begin{center}
A rational double point has finite local fundamental group.
\end{center}

For the following see \cite{B3} Satz 2.8, \cite{Prill} Thm. 3.
From a merely topological point of view, $x \in X$ possesses a
neighbourhood $U$ with $U^{\prime} := U \setminus \{ x \}$ having
a finite universal cover $V^{\prime} \twoheadrightarrow
U^{\prime}$. This can be uniquely extended to a ramificated cover
$V \twoheadrightarrow U$ by adding a point $y$ to $V^{\prime}$.
Moreover, we can equip $V$ with a normal analytic structure such
that $ V \twoheadrightarrow U$ becomes an analytically ramificated
cover. Since $V^{\prime}$ is simply
connected, we see $\pi _{V,y}=1$. By another fundamental theorem
in Mumford's paper (\cite{Mf} p. 18), this shows the non-singularity
of $V$ at $y$. Now $\pi _{X,x}$ is operating via cover
transformations on $V^{\prime}$, hence also on $V$ with fixed
point $y$. We need another definition to state our results so far.

\begin{defprop}{\rm  (Two-dimensional quotient singularities)}
For a neighbourhood $V$ of the origin $O$ in $\CC^2$ and a finite group $G$ of analytic automorphisms of $V$ fixing $O$, the quotient space
$V / G$ has the structure of a normal complex analytic surface and the projection $V \twoheadrightarrow V / G$ is analytic \cite{B3}.
We say that $U$ is a {\em two-dimensional quotient singularity},  if $U$ is isomorphic to a singularity of the form $V/G$.
\end{defprop}

We have just seen:

\begin{center}
A rational double point is a quotient singularity.
\end{center}

By a simple linearization argument (\cite{B3} Lemma 2.2), we can restrict ourselves to the study of quotient singularities of the form $\CC^2 / G$ where
$G$ is a finite subgroup of $\text{GL}(2,\CC)$. Obviously, conjugated subgroups yield isomorphic quotient spaces.
The conjugacy classes of finite subgroups of $\text{GL}(2,\CC)$ have been listed by Du Val (\cite{DuVal3} \S 21).
Later, Prill showed that only a particular class of finite subgroups has to be studied: the so-called {\em small} subgroups \cite{Prill}.
He also classified them (\cite{Prill} Satz 2.3). Using Prill's results, Brieskorn gave a complete classification of the quotient spaces that can
arise in terms of the Dynkin diagram of their resolution (\cite{B3} p. 348). This shows that two-dimensional quotient singularities are taut and finishes
Brieskorn's proof.

% ----------------------------------------------------------------------------------------------------------------------------------------------------------------------------------------------------------------------------------------
\section{Seven characterizations of rational double points}

The following remark allows us to use the intermediate results in
the above discussion \ref{bries} to give alternative
characterizations of rational double points in the analytic
category.

\begin{remark}
\label{lastbit}
It can be shown that $\CC^2 / G$, for $G \subset \text{\rm GL}(2,\CC)$ finite,
embeds in codimension one if and only if $G$ is a subgroup
of $\text{\rm SL}(2,\CC)$ {\rm (\cite{Df} cor. 5.3)}.
\end{remark}

\begin{theorem}
\label{characterizations}
{\rm (\cite{Df})}

Let $(X,x)$ be a normal surface singularity that {\bf embeds in codimension one}. Then the following conditions are equivalent in the
complex analytic category.
\begin{enumerate}
\item $(X,x)$ is a rational double point.
\item $(X,x)$ has a good resolution with an exceptional set consisting of rational curves with self-intersection-number $-2$.
\item $(X,x)$ has a good resolution with an exceptional set consisting of rational curves and a Dynkin diagram listed in table \ref{rdptable}.
\item The local fundamental group of $(X,x)$ is finite.
\item $(X,x)$ is a two-dimensional quotient singularity.
\item $(X,x)$ is isomorphic to $\CC^2 / G$ for finite $G \subset \text{\rm GL}(2,\CC)$.
\item $(X,x)$ is isomorphic to one of the affine varieties studied in section \ref{example1}.
\end{enumerate}

\end{theorem}

{\bf Remark on the proof:} We have already seen a proof of the
implications $(1) \Rightarrow (2) \Rightarrow (3) $ and presented
some ideas for $ (3) \Rightarrow (4) \Rightarrow (5) \Rightarrow
(6)$. In the example \ref{icosahedral}, we studied a special case
of $(7) \Rightarrow (1)$. The last implication $(6) \Rightarrow
(7)$ is precisely remark \ref{lastbit}.

\vspace{3mm} In his survey article \cite{Df}, Durfee presents
further ten characterizations of rational double points. The
characterizations he gives provide a connection of rational double
points for example with weighted homogeneous
polynomials, vanishing cycles, a certain limit
involving volumes, monodromy groups and Morse functions.
A more number-theoretical characterization in terms
of almost factorial rings (fast-faktoriell) rings is due to
Brieskorn (\cite{B3} Satz 1.5).
Finally, a link with elementary catastrophes is discussed in a survey article by Slodowy (\cite{Sl} 9).

% -----------------------------------------------------------------------------------------------------------------------------------------------------------------------------------------------------------------------------------------------
\section{Example: The conical double point}

Although we did not give a proof of theorem
\ref{characterizations}, we shall at least study an example to
illustrate the phenomena encountered there. We will consider the
{\em conical double point} $X:=V(f) \subset \CC^3$ where $f=xz -
y^2 \in \CC [x,y,z]$. Note that after a change of variables $xz -
y^2$ becomes $x_1^2 + x_2^2 +x_3^2$; so $X$ is just the surface
singularity labeled $A_1$ from table \ref{rdptable}. $X$ is a
double cone with vertex a rational double point:

\begin{center}
\cone.
\end{center}

Obviously, $X$ has a normal singularity at $x=(0,0,0)$ and $X$ is embedded in codimension one.
$(X,x)$ is a double point, since

$$\text{rank}_{\CC}\frac{(x,y,z)^k + (f)}{(x,y,z)^{k+1} +(f)} = 2k+1,  $$

i.e. the leading coefficient of the Hilbert-Samuel polynomial of the local ring $\OO _{X,x}$ at  $x$ it two.

Furthermore, $(X,x)$ is absolutely isolated, because the singularity can be resolved by a single blow-up at $x$, as can be easily seen by the
toric description of $X$ \cite{Fulton}:

\begin{center}
\toric.
\end{center}

We can work out the blow-up explicitly and obtain
\label{blowup}

$$ \widetilde{X}^{\prime} \subset \mathop{{\mathbb A}^3 \times {\mathbb P}^2}_{(x,y,z;p:q:r)} $$

given by equations

$$ xz - y^2, pr - q^2, py=qx, pz=rx, qz=ry .$$

$\widetilde{X}^{\prime}$ is isomorphic to

$$ \widetilde{X} \subset \mathop{{\mathbb A}^2 \times {\mathbb P}^1}_{(x,z;u:v)} $$

cut out by

$$ xv^2 = zu^2  $$

via

\begin{eqnarray*}
\widetilde{X} & \stackrel{\sim}{\longrightarrow} & \widetilde{X}^{\prime} \\
(x,z;u:v) & \mapsto & (x,\frac{v}{u}x\text{ or }\frac{u}{v}z,z;u^2:uv:v^2)
\end{eqnarray*}

Now $\widetilde{X}$ is a line bundle on \PP\, whose zero-section has self-intersection-number $-2$. Thus $\widetilde{X}$  is just the line bundle
on the projective line associated to the sheaf $\OO _{\PP}(-2)$, i.e. the cotangent bundle $T^{\ast} \PP$.
See also \cite{Lm} IV \S7.

Since $\widetilde{X}$ is a smooth variety, we have a resolution of $(X,x)$

$$ \pi : \widetilde{X} \cong \widetilde{X}^{\prime} \twoheadrightarrow X . $$

The exceptional set $E$ of  $\pi$ is precisely the zero-section of  $\widetilde{X}$, hence isomorphic to \PP. Moreover $E^2 = -2 $.
The real picture reflects the situation very nicely

\begin{center}
\resreal.
\end{center}

For the numerical cycle we get $\Znum = E$, i.e. $(X,x)$ is
rational by theorem \ref{theo2}, and characterization $(1)$ is
verified.

As already mentioned in section \ref{C2G}, $X$ is isomorphic to the affine orbit variety $\CC^2 / \{ \pm 1 \}$, where we write $-1$ for the reflection in the
origin of the complex plane $\CC^2$. This corresponds to characterizations $(5)$, $(6)$ and $(7)$.

Finally, let us calculate the local fundamental group $\pi_{X,x}$. We have a covering map

$$ \CC^2 \setminus \{O\} \twoheadrightarrow X \setminus \{ x\} $$

with covering transformation group $\{ \pm 1 \}$.
For every $\{ \pm 1 \}$-invariant simply connected neighbourhood
$U$ of $O \in \CC^2$, we observe that $U \setminus \{ O \}$ is
also simply connected, hence
$$
\pi_1 ((U \setminus \{ O \}) / \{ \pm 1 \}) = \{ \pm 1 \} .
$$
But such an $U$ can be chosen arbitrarily small, thus
$$
\pi_{X,x} = \{ \pm 1\},
$$
which is finite and shows characterization $(4)$.

% -------------------------------------------------------------------------------------------------------------------------------------------------------------------------------------------------------------------------------------------------
\section{Lie groups and rational double points}

The "$A_n$-$D_n$-$E_n$" - labeling of the various types of
rational double points was actually borrowed from the
classification theory of Lie groups. In this last section we will
sketch some of the deep connections between Lie groups and
rational double points. Essentially, we shall
give a summary of \cite{Sl} 10.

% -------------------------------------------------------------------------------------------------------------------------------------------------------------------------------------------------------------------------------------------------
\subsection{Dynkin diagrams of simple Lie groups}

A connected complex Lie group is called {\it (almost) simple}, if it
contains no normal subgroup of positive dimension.
In their classification theory, the simply connected simple Lie
groups play a special r\^{o}le as their universal coverings
(which are finite (\cite{Sl} 10)).
These groups are classified by their corresponding
Dynkin diagrams \cite{FH}. Surprisingly, the Dynkin diagrams of
table \ref{rdptable} occur again.

We recall briefly the relevant part of this classification; note that most of the
following facts hold in a more general context (\cite{FH},
\cite{Sl2} 3.1).

Let $G$ be a simply connected simple Lie group of rank $r$ and
\GG\ its Lie algebra. We fix a maximal torus $T\cong(\CC^{\ast})^r$
of $G$ with {\it character group}

$$ X^{\ast} ( T ) = \text{Hom} ( T, \CC^{\ast}) = \ZZ^r . $$

We denote the normalizer of $T$ in $G$ by $N_G(T)$. The group
$W:=N_G(T) / T$ is finite and is called the {\it Weyl group} of $G$
with respect to $T$.

The restriction of the adjoint representation of $G$ on \GG\ to $T$
has eigenspaces $\GG_{\alpha}$ on which $T$ acts by the character
$\alpha \in X^{\ast}(T)$ and we obtain the {\it Cartan
decomposition} of \GG\

$$ \GG = \bigoplus_{\alpha \in X^{\ast}(T)} \GG_{\alpha} . $$

The finite set $\Sigma := \{ \alpha \in X^{\ast}(T) : \alpha \neq
0 , \GG_{\alpha} \neq \{0\}\} $ is called the {\it root space}.
Clearly, $\Sigma$ is invariant under the action of $W$.

We can define a $W$-invariant scalar product $\langle\cdot,\cdot\rangle$ on $X^{\ast}(T)$
(called the {\it Killing form}) such that the elements of $W$
become reflections in the hyperplane perpendicular to a root
$\alpha$

$$ \beta \mapsto \beta - \frac{2 \langle \alpha, \beta \rangle}{\langle
\alpha , \alpha \rangle} \alpha . $$

Here $\frac{2 \langle \alpha, \beta \rangle}{\langle
\alpha , \alpha \rangle}$ must be integer.
(The main tool in proving this  and similar facts is identifying
the subalgebra $\GG_{\alpha} \oplus \GG_{- \alpha} \oplus [ \GG_{\alpha},\GG_{-\alpha} ]$
with $\mathfrak{sl}(2,\CC)$ and applying the representation
theory of $\mathfrak{sl}(2,\CC)$.)

The geometry of how $\Sigma$ sits in the Euclidean lattice $(X^{\ast}(T),\langle \cdot,
\cdot \rangle)$ is very rigid. For example

$$ 4 \cos^2 \angle (\alpha, \beta) =\frac{ 4 \langle \alpha, \beta
\rangle^2 }{\langle \alpha, \alpha \rangle \langle \beta, \beta
\rangle}$$

must be an integer between zero and four, i.e. there are just a
few possibilities for the angle between two roots $\alpha$ and
$\beta$.

By choosing a direction $l \in T$ (in general position with
respect to $\Sigma$) we can specify the {\it positive roots}
$\alpha \in \Sigma$ to be those with $\alpha(l) > 0$.
In particular, we can focus on {\it simple roots}: these are
positive roots that are not the sum of two other positive roots.
The system of simple roots gives rise to a Dynkin diagram, where
we take a vertex for each simple root and join two vertices by
exactly
$ 4 \cos^2 \angle (\alpha, \beta) $
lines. If we insist on all roots $\alpha$ having the same
length $\langle \alpha, \alpha \rangle$, the only possibilities for
the Dynkin diagram are

\begin{center}
$A_n$ \raisebox{-2mm}{\Anplain} (n vertices), \\
$D_n$ \raisebox{-7mm}{\Dnplain} (n vertices), \\
$E_6$ \raisebox{-7mm}{\Esechsplain}, \\
$E_7$ \raisebox{-7mm}{\Esiebenplain} and\\
$E_8$ \raisebox{-7mm}{\Eachtplain}.
\end{center}

These diagrams $A_n$, $D_n$, $E_6$, $E_7$ and $E_8$ actually occur
for the simply connected simple Lie groups corresponding to the
classical Lie algebras $\mathfrak{sl}(n+1,\CC)$,
$\mathfrak{so}(2n,\CC)$ and the exceptional Lie algebras
 $\mathfrak{e}_6$, $\mathfrak{e}_7$ and $\mathfrak{e}_8$, respectively.

It can be shown that a simply connected simple Lie group can be recovered from
its Dynkin diagram.

% -------------------------------------------------------------------------------------------------------------------------------------------------------------------------------------------------------------------------------------------------
\subsection{A theorem of Brieskorn}

Let $G$ be a simply connected simple Lie group.
We consider the quotient $H$ of $G$ by its adjoint action in the
category of algebraic varieties $p: G \twoheadrightarrow H$. There
is an explicit way to describe $p$. Let $r$ be the rank of $G$ and $\rho_i: G \rightarrow
\text{GL}(V_i), i=1, \dots, r$ be the $r$ fundamental irreducible
representations of $G$ on finite-dimensional vector spaces. Then the character map

\begin{eqnarray*}
\chi: G & \rightarrow & \CC^r \\
g & \mapsto & (\dots,\text{trace}_{V_i} \rho_i(g),\dots)
\end{eqnarray*}

coincides with $p$.

\paragraph{Example:} For $G=\text{SL}(n,\CC)$ we have $\text{rank }G=n-1$ and the $n-1$
fundamental irreducible representations are given by the exterior powers

$$V_i = {\bigwedge}^i\, \CC^n . $$

The corresponding characters are, up to sign, just the
non-trivial coefficients of the characteristic polynomial

$$\text{char}(g)=\text{det}(\lambda - g) = \lambda^n -
\text{trace}(g) \lambda ^{n-1}
+ \text{trace}({\wedge}^2 g)\lambda^{n-2} -+ \dots .$$

Thus we can regard $\chi$ as associating to $g \in
\text{SL}(2,\CC)$ its characteristic polynomial.

\vspace{3mm}

The point is now to study the {\it unipotent variety}

$$\text{Uni}(G) := p^{-1}(p(e)). $$

\paragraph{Example:} For $G=\text{SL}(n,\CC)$ the unipotent
variety consists precisely of the unipotent matrices.

\vspace{3mm}

The variety $\text{Uni}(G)$ is a finite union of conjugacy classes
and contains a unique conjugacy class of dimension $d:=\text{dim }G
-r$ (since $p$ is flat) --- the {\it regular class}.
The complement of the regular class in $\text{Uni}(G)$ is the
closure of a unique conjugacy class of dimension $d-2$ --- the
{\it subregular class} $\text{Sub}(G)$.

We choose a $S \subseteq G$ such that

\begin{eqnarray*}
& S \text{ is smooth of dimension }\text{dim }G-d+2, \\
& S \cap \text{Sub}(G)=\{ x \} \text{ and } \\
& T_x S + T_x \text{Sub}(G) = T_x G,
\end{eqnarray*}

i.e. we require that $S$ is a {\it slice} of codimension $d-2$ transversal to
$\text{Sub}(G)$ at an element $x \in \text{Sub}(G)$.

Let $X:=S \cap \text{Uni}(G)$.

The following theorem was conjectured by Grothendieck and proved
by Brieskorn.

\begin{theorem}
\label{lie}
{\rm (\cite{B5}, \cite{Sl} 10)}
If $G$ is a simply connected simple Lie group of type $A_n$, $D_n$
or $E_n$, then $(X,x)$ is a rational double point.
\end{theorem}

\paragraph{Example:} We want to illustrate this theorem in the
simplest possible case $G=\text{SL}(2,\CC)$. All regular unipotent
elements are conjugate to the matrix
$$
\left(
\begin{array}{cc}
1 & 1 \\
0 & 1
\end{array}
\right)
$$
and there exists just a single subregular unipotent element
$$
x:=
\left(
\begin{array}{cc}
1 & 0 \\
0 & 1
\end{array}
\right) .
$$
As a transversal slice we can simply take $S:=\text{SL}(2,\CC)$.
Now $p=\chi$ is given by the trace
\begin{eqnarray*}
\chi : \text{SL}(2,\CC) & \rightarrow & \CC \\
\left(
\begin{array}{cc}
a & b \\
c & d
\end{array}
\right)
& \mapsto & a+d
\end{eqnarray*}
and we get
\begin{eqnarray*}
X & = & S \cap \text{Uni}(G) \\
& =& \text{Uni}(G)=\chi^{-1}
\left(
\chi
\left(
\begin{array}{cc}
1 & 0 \\
0 & 1
\end{array}
\right)
\right) \\
& = & \left\{
\left(
\begin{array}{cc}
1+x & y \\
z & 1+u
\end{array}
\right) : x+u=0 \text{ and } xu-yz=0
\right\} \\
& = & \{ (x,y,z) : x^2 + yz = 0\},
\end{eqnarray*}

i.e. $X$ has a conical double point at
$x=(0,0,0)$.

% -------------------------------------------------------------------------------------------------------------------
\subsection{Resolutions of rational double points in the Lie group context}

A closed subgroup $P\subseteq G$ is called {\it parabolic}, if the
quotient space $G / P$ is a projective variety. The minimal parabolic
subgroups are the {\it Borel subgroups}. All Borel subgroups are
conjugate to each other in $G$ and the normalizer $N_G(P)$ of a
parabolic subgroup coincides with P. Thus the set of all Borel
subgroups $\mathcal B$ becomes a projective variety

$$ {\mathcal B} = G / B $$

where $B$ is any Borel subgroup of $G$. More generally, ${\mathcal
P}:= G/P$ may be identified with the set of subgroups conjugate to
the parabolic subgroup $P$.

\paragraph{Example:} The parabolic subgroups of $\text{SL}(n,\CC)$
are exactly the stabilizer of the flags

$$ 0 \subset V_{i_1} \subset \dots \subset V_{i_k} \subset \CC^n
\text{ with } \text{rank}_{\CC}V_{i_j}=i_j, j=1,\dots,k .$$

Hence the Borel subgroups correspond to maximal flags

$$ 0 \subset V_1 \subset \dots \subset V_{n-1} \subset \CC^n , $$

i.e. are conjugate to the subgroup of upper triangular matrices.

\vspace{3mm}

It was Springer who showed that the natural projection from the
incidence variety

$$ I := \{ (x,B) \in \text{Uni}(G) \times {\mathcal B} : x \in B
\} $$

to $\text{Uni}(G)$ is a $G$-equivariant resolution of the
singularities of $\text{Uni}(G)$

$$ \pi : I \twoheadrightarrow \text{Uni}(G) . $$

Let $G,x,S,X$ be as in theorem \ref{lie}. It is a consequence of
the $G$-equivariance of $\pi$ that the restriction

$$ \pi : \widetilde{X} := \pi^{-1}(X) \twoheadrightarrow X $$

is again a resolution. In fact, it is a {\it minimal} one, that
is, we cannot apply theorem \ref{app4} to obtain a resolution with smaller
exceptional set.

We can interpret the exceptional set

$$ E:= \pi^{-1} (x) $$

in two different ways.

On one hand, we know from theorem \ref{lie}, that $(X,x)$ is a
rational double point. Hence $E$ must be a bunch of projective
lines \PP\ intersecting each other as prescribed by the Dynkin
diagram $\Gamma$ of $(X,x)$.

On the other hand, we can write

$$ E = \{ (x,B) \in \{x\} \times {\mathcal B} : x \in B\} . $$

The vertices of the Dynkin diagram $\Gamma_G$ of $G$ correspond to
the simple roots of $G$ (after a maximal torus $T_0$ and a
direction $l \in T_0$, or equivalently, a Borel subgroup $B_0
\supset T_0$ have been specified). Let $P_{\alpha}$ be the minimal
proper (i.e. non-Borel) parabolic subgroup generated by $B_0$ and
the root subgroup $U_{-\alpha}$, where $\alpha$ is a simple root.
Because of $N_G(P_{\alpha})=P_{\alpha}$, we can identify the set
of subgroups conjugate to $P_{\alpha}$ with the projective
variety

$$ {\mathcal P}_{\alpha}:= G / P_{\alpha} .$$

The natural map

$$ f_{\alpha} : {\mathcal B} \cong G/B_0 \rightarrow {\mathcal
P}_{\alpha} \cong G/P_{\alpha} $$

maps each Borel subgroup $B$ to the unique parabolic subgroup $P\in
{\mathcal P}_{\alpha}$ containing $B$. Since each $P \in {\mathcal
P}_{\alpha}$ has semisimple rank 1, this map has projective lines
as fibres.

Steinberg and Tits showed that $E$ is a bunch of projective lines
--- one line of the form $f_{\alpha}^{-1}(P), P\in {\mathcal
P}_{\alpha}$ for every simple root $\alpha$ --- which intersect as
prescribed by the edges of $\Gamma_G$.

\paragraph{Example:} We verify these statements by explicit
calculation in the simplest non-trivial case
$G=\text{SL}(3,\CC)$.

As maximal torus $T_0$ we may take the diagonal matrices and as
Borel subgroup the upper triangular matrices.

The root space of $\text{SL}(3,\CC)$ is
\begin{center}
\roots
\end{center}
where the character $\alpha_{i,j} \in X^{\ast}(T)$ is defined by

$$ \alpha_{i,j}
\left(
\begin{array}{ccc}
a_1 && \\
& a_2 & \\
&& a_3
\end{array}
\right)
=
a_i - a_j . $$

The eigenspace $\mathfrak{sl}(3,\CC)_{\alpha_{i,j}}$, for
example, consists of those matrices
$(a_{i,j})_{i,j=1,\dots,3}\in\mathfrak{sl}(3,\CC)$
whose single non-zero entry is $a_{2,1}$. Using the
exponential map
$$
\text{exp} : \mathfrak{sl}(3,\CC) \rightarrow \text{SL}(3,\CC)
$$
we see that the root subgroup $U_{\alpha_{21}}$ is
$$
U_{\alpha_{21}} =
\left\{
\left(
\begin{array}{ccc}
1 & 0 & 0 \\
\lambda & 1 & 0 \\
0 & 0 & 1
\end{array}
\right) : \lambda \in \CC
\right\} .
$$

The simple roots are $\alpha_{1,2}$ and $\alpha_{2,3}$ and we
obtain
\begin{center}
\Azweibig\
\end{center}
as Dynkin diagram of $\text{SL}(3,\CC)$ (\cite{FH} 12).
The unipotent variety $\text{Uni}(\text{SL}(3,\CC))$ is given by
the unipotent matrices, all regular unipotent elements are
conjugate to
$$
\left(
\begin{array}{ccc}
1 & 1 & 0 \\
0 & 1 & 1 \\
0 & 0 & 1
\end{array}
\right)
$$
and all subregular unipotent elements are conjugate to
$$
x:=\left(
\begin{array}{ccc}
1 & 0 & 0 \\
0 & 1 & 1 \\
0 & 0 & 1
\end{array}
\right) .
$$
We have the minimal proper parabolic subgroups
$$
P_{\alpha_{1,2}} = \langle B_0 , U_{\alpha_{2,1}} \rangle
$$
and
$$
P_{\alpha_{2,3}} = \langle B_0 , U_{\alpha_{3,2}} \rangle
$$
which are the stabilizer of the flags
$$
0 \subset \text{span}_{\CC}\{ v_1 , v_2 \} \subset \text{span}_{\CC}\{ v_1 , v_2 ,v_3
\}
$$
and
$$
0 \subset \text{span}_{\CC}\{ v_1 \} \subset \text{span}_{\CC}\{ v_1 , v_2 ,v_3
\},
$$
respectively. Clearly, $P_{\alpha_{2,3}}$ is conjugate to the
parabolic subgroup $P_{\alpha_{2,3}}^{\prime}$ stabilizing the flag
$$
0 \subset \text{span}_{\CC}\{ v_2 \} \subset \text{span}_{\CC}\{ v_1 , v_2 ,v_3
\}.
$$
The set of Borel subgroups containing $x$, which we had identified
with $E$, is given by
$$
f_{\alpha_{1,2}}^{-1}(P_{\alpha_{1,2}}) \cup
f_{\alpha_{2,3}}^{-1}(P_{\alpha_{2,3}}^{\prime}).
$$
(All the other fibres of the maps $f_{\alpha_{1,2}}$ and
$f_{\alpha_{2,3}}$ contain only a finite number of Borel subgroups
which contain $x$.)

The two fibres intersect in a single Borel
subgroup: the subgroup stabilizing the maximal flag
$$ 0 \subset \text{span}_{\CC}\{ v_2  \} \subset \text{span}_{\CC}\{ v_1 , v_2 \} \subset \text{span}_{\CC}\{ v_1 , v_2 , v_3
\}  . $$
Hence the Dynkin diagram of $E$ looks like
\begin{center}
\Azweibig .
\end{center}

% -------------------------------------------------------------------------------------------------------------------------------------------------------------------------------------------------------------------------------------------------
\appendix
\section{Appendix - results from Algebraic Geometry}

\begin{theorem}
\label{app1}
{\rm (Riemann-Roch on a surface)(\cite{Ha} V.1.6)}
If $D$ is any divisor on the non-singular surface $X$, then

$$ \chi(\Ox(D)) = \frac{1}{2} D \cdot (D-K) + 1 + p_a(X) $$

where $K$ is a {\em canonical divisor} on $X$ (\cite{Ha} V.1.4.4).
\end{theorem}

A simple corollary is the following.

\begin{theorem}
\label{app2}
{\rm (General adjunction formula)(\cite{Ha} Ex. V.1.3.(a))}
If $D$ is an ef\-fec\-tive divisor on the non-singular surface $X$, then

$$ 2 p_a(D) -2 = D \cdot (D+K) . $$
\end{theorem}

\begin{theorem}
\label{app3}
{\rm (Grothendieck's theorem on formal functions)(\cite{EGA3} 4.2.1 or for projective morphisms \cite{Ha} III.11.1)}
Let $f: X \rightarrow Y$ be a proper morphism of noetherian schemes and $\mathcal F$ a coherent sheaf on $X$. For $y \in Y$ denote by
$m_y \subset \OO _{Y,y}$ the maximal ideal of the stalk at $y$. Then we have a natural isomorphism

$$ 0 = \left( \left( R^i f_{\ast} {\mathcal F} \right)_y \right)\widehat{}\ \cong\
  \mathop{\varprojlim}_{k=1}^{\infty} H^i\left( f^{-1}(y), {\mathcal F} \otimes _{\OO _Y} {\OO_{Y,y} \over {m_y}^k} \right)
\text{ for all } i \in {\mathbb N}$$

where the completion is taken with respect to the $m_y$-adic topology.
\end{theorem}

\begin{theorem}
\label{app4}
{\rm (Castelnuovo's criterion for contracting a curve)(\cite{Lm} IV \S 15, \cite{Ha} V.5.7)}
If $C$ is a curve on a non-singular surface $X$ with $C \cong \PP$ and $C^2=-1$, then
there exists a morphism $f: X \rightarrow X^{\prime}$ to a non-singular surface $X^{\prime}$ which contracts $C$ to a point $p$,
such that $X$ is isomorphic via $f$ to the
blow-up of $X^{\prime}$ with center $p$, and $C$ is the exceptional curve.
\end{theorem}

\begin{theorem}
\label{app5} {\rm (Zariski's connectedness theorem)(\cite{EGA3}
4.3.1, \cite{Ha} III.11.4)} Let $f: X \rightarrow Y$ be a
birational morphism between projective varieties and assume that
$Y$ is normal. Then $f$ has connected fibers.
\end{theorem}

\begin{theorem}
\label{app6}
{\rm (A vanishing theorem of Grothendieck)(\cite{Ha} III.2.7)}
For any sheaf of abelian groups ${\mathcal F}$ on a noetherian scheme $X$ of dimension $n$, we have

$$ H^i(X,{\mathcal F}) =0 \text{ for } i>n . $$
\end{theorem}

\begin{theorem}
\label{app7} {\rm (A vanishing theorem for higher direct image
sheaves)(\cite{Ha} III.11.2)} Let $f: X \rightarrow Y$ be a
projective morphism of noetherian schemes and denote by $r$ the
maximal dimension of its fibers. Then for all coherent sheaves
${\mathcal F}$ on $X$, we have

$$ R^i f_{\ast}\mathcal F = 0 \text{ for } i>r .$$
\end{theorem}

\newpage

% -----------------------------------------------------------------------------------------------------------------------

\end{document}